\documentclass[letterpaper,11pt]{article}

\usepackage{amsthm}
\usepackage{amsmath}

\usepackage{mathrsfs}
\newcommand{\fancyfont}[0]{\mathscr}

\newtheorem{thm}{Theorem}[section]

\newtheorem{lemma}[thm]{Lemma}
\newtheorem{prop}[thm]{Proposition}
\newtheorem{cor}[thm]{Corollary}

\newtheorem{conj}[thm]{Conjecture}

\newcommand{\beq}[1]{\begin{equation}\label{#1}}
\newcommand{\enq}[0]{\end{equation}}

\newcommand{\bn}[0]{\bigskip\noindent}
\newcommand{\mn}[0]{\medskip\noindent}
\newcommand{\nin}[0]{\noindent}

\newcommand{\sub}[0]{\subseteq}
\newcommand{\sm}[0]{\setminus}
\renewcommand{\dots}[0]{,\ldots,}

\newcommand{\ov}[0]{\overline}

\newcommand{\cee}[0]{{\fancyfont C}}

\newcommand{\f}[0]{{\cal F}}
\newcommand{\g}[0]{{\cal G}}
\newcommand{\h}[0]{{\cal H}}

\newcommand{\I}[0]{{\cal I}}

\newcommand{\J}[0]{{\cal J}}
\newcommand{\K}[0]{{\cal K}}
\newcommand{\m}[0]{{\cal M}}
\newcommand{\N}[0]{{\cal N}}

\newcommand{\pee}[0]{{\cal P}}
\newcommand{\Q}[0]{{\cal Q}}
\newcommand{\R}[0]{{\cal R}}

\newcommand{\T}[0]{{\cal T}}

\newcommand{\W}[0]{{\cal W}}

\newcommand{\ra}[0]{\rightarrow}

\newcommand{\ZZ}[0]{{\bf Z}}

\newcommand{\Nn}[0]{{\bf N}}

\newcommand{\XX}[0]{{\bf X}}

\newcommand{\Tr}[0]{{\rm Tr}}

\newcommand{\cc}[0]{\mbox{{\sf c}}}
\newcommand{\ddd}[0]{\mbox{{\sf d}}}
\newcommand{\ttt}[0]{\mbox{{\sf t}}}
\newcommand{\aaa}[0]{\mbox{{\sf a}}}
\newcommand{\bbb}[0]{\mbox{{\sf b}}}
\newcommand{\ww}{\mbox{{\sf w}}}

\newcommand{\0}[0]{\emptyset}

\renewcommand{\qed}[0]{\begin{flushright} \rule{2mm}{3mm} \end{flushright}}
\def\qqed{\null\nobreak\hfill\hbox{${\diamondsuit}$}\par\smallskip}
\def\qqqed{\null\nobreak\hfill\hbox{\rule{2mm}{3mm} }\par\smallskip}

\newcommand{\C}[0]{\binom}

\newcommand{\Cc}[0]{\tbinom}
\newcommand{\ga}[0]{\alpha }
\newcommand{\gb}[0]{\beta }
\newcommand{\gc}[0]{\gamma }
\newcommand{\gd}[0]{\delta }
\newcommand{\gD}[0]{\Delta }

\newcommand{\gl}[0]{\lambda }

\newcommand{\gs}[0]{\sigma}

\newcommand{\gz}[0]{\zeta}
\newcommand{\eps}[0]{\varepsilon }
\newcommand{\vt}[0]{\vartheta}

\newcommand{\vr}[0]{\varrho}

\newcommand{\ii}[0]{\I^*}
\newcommand{\lct}[0]{l^{\C{t}{2}}}
\newcommand{\ct}[0]{\C{t}{3}}
\newcommand{\nsa}[0]{\not\sim_A}
\newcommand{\abc}[0]{(A,B,C)}
\newcommand{\Cl}[0]{{\cal K}}
\newcommand{\Clp}[0]{{\cal K}(\pee^*)}
\newcommand{\cone}[0]{c_1}
\newcommand{\ctwo}[0]{c_2}
\newcommand{\hhh}[0]{r}

\newcommand{\es}[0]{{2\eps_2 l}}

\newcommand{\sugg}[1]{}

\begin{document}

\renewcommand{\thefootnote}{\fnsymbol{footnote}}
\footnotetext{AMS 2010 subject classification:  05A16, 05C65, 68R05}
\footnotetext{Key words and phrases:  k-SAT function,
asymptotic enumeration, hypergraph regularity lemma
}
\title{The number of 3-SAT functions\footnotemark }

\author{
L. Ilinca and J. Kahn}
\date{}

\footnotetext{ * Supported by NSF grant DMS0701175.}

\maketitle

\begin{abstract}
With $G_k(n)$ the number of functions
of $n$ boolean variables definable by
$k$-SAT formulae,
we prove that $G_3(n)$ is asymptotic to $2^{n+\C{n}{3}}$.
This is a strong form of the case $k=3$ of a conjecture of Bollob\'as, Brightwell and
Leader stating that for fixed $k$, $\log_2G_k(n)\sim \C{n}{k}$.
\end{abstract}

\section{Introduction}\label{Intro}

Let $X_n=\{x_1,\ldots,x_n\}$ be a collection of Boolean variables.
Each variable $x$ is associated with a
${\it positive}$ literal, $x$, and a ${\it negative}$ literal $\bar{x}$.
Recall that a {\it k-SAT formula} (in disjunctive normal form)
is an expression $\cee$ of the form
\beq{kS}
C_1\vee \cdots \vee C_t,
\enq
with $t$ a positive integer and each $C_i$ a {\em k-clause};
\glossary{name={$k$-SAT formula},sort=SAT}
that is, an expression
$y_1\wedge \cdots \wedge y_k$, with $y_1\dots y_k$ literals corresponding to
different variables.
\glossary{name={$k$-clause},sort=clause}
A formula (\ref{kS}) defines a Boolean function of $x_1\dots x_n$
in the obvious way;
any such function is a {\em k-SAT function}.
Though we will be concerned here almost exclusively with the case $k=3$,
we leave the discussion general for the moment.

Following \cite{BBL}, we write $G_k(n)$ for
the number of $k$-SAT functions of $n$ variables.
\glossary{name={$G_k(n)$},description={The number of $k$-SAT functions of $n$ variables.},sort=G}
Of course $G_k(n)$ is at most
$\exp_2[2^k\C{n}{k}]$,
the number of $k$-SAT formulas; on the other hand it's
easy to see that
\beq{3S}
G_k(n) > 2^n(2^{\C{n}{k}}-n2^{\C{n-1}{k}}) \sim 2^{n+\C{n}{k}}
\enq
(all formulas obtained by choosing
$y_i\in \{x_i,\bar{x}_i\}$ for each $i$
and a set of clauses using precisely the literals
$y_1\dots y_n$ give different
functions).

The problem of estimating $G_3(n)$ was suggested by
Bollob\'as, Brightwell and Leader \cite{BBL}.
They showed
\beq{bbl}
G_k(n)\leq \exp_2[(2 \sqrt{\pi})\Cc{n}{k}],
\enq
for $k< n/2$
and conjectured that
\beq{bblconj}
\log_2G_k(n)<(1+o(1))\Cc{n}{k}.
\enq

\nin
for any fixed $k$.
Even $k=2$ is not easy; here (\ref{bblconj})
was proved in \cite{BBL}, and the precise asymptotics---
\beq{allen}
G_2(n)\sim \exp_2[n+\Cc{n}{2}]
\enq
---conjectured in \cite{BBL}
were proved in \cite{Allen} and (later) in \cite{IK}.
As is often the case, nothing from this earlier work
seems to be of much help in treating larger $k$.

Here, for $k=3$, we
prove (\ref{bblconj}) and more, again showing (as in \eqref{allen})
that (\ref{3S}) gives the
asymptotics not just of $\log G_3(n)$, but of $G_3(n)$ itself:
\begin{thm}
\label{SATthm}
$~~~~G_3(n)
\sim 2^{n+\C{n}{3}}.$
\end{thm}

For a formula $\cee$ as in (\ref{kS}) we may identify the
associated function, say $f_{\cee}$, with the set
(henceforth also referred to as a ``$k$-SAT function")
$F(\cee)\sub \{0,1\}^n$ of satisfying assignments for $\cee$
(that is, $F(\cee) = f_{\cee}^{-1}(1)$).
For our purposes it will also usually be convenient to think
of $\cee$ as the set $\{C_1\dots C_t\}$ of clauses.
Then $F(\cee')\sub F(\cee)$ whenever $\cee'\sub \cee$, and
we say $\cee$ is {\em irredundant} if it is a minimal formula
giving $F(\cee)$; that is, if
$F(\cee')\subset F(\cee)$ for each $\cee'\subset \cee$.
\glossary{name={irredundant}}
Of course each 3-SAT function $F$ corresponds to at least one
irredundant $\cee$, so that,
with $I(n)=I_3(n)$ denoting the number of irredundant formulas on $X_n$,
\glossary{name={$I(n)$},description={$=I_3(n)$; see $I_k(n)$.},sort=I}
Theorem \ref{SATthm} is contained in
\begin{thm}\label{Irr}
$~~I(n)
\sim 2^{n+\C{n}{3}}$.
\end{thm}

This (together with (\ref{3S})) says that in fact most
$F$'s admit only one irredundant formula.
We regard this simple idea as one of the keys to the present work:
it allows us to forget about functions and work directly
with formulas, which are easier (though to date still not easy) to handle.

Notice that $\cee$ is irredundant iff for each $C\in \cee$ there is
some (not necessarily unique) {\em witness}
$\ww_C\in \{0,1\}^n$ that satisfies $C$ but no other clause in $\cee$
(i.e. $\ww_C\in F(\cee)\sm F(\cee\sm\{C\})$).
\glossary{name={witness}}
Such witnesses will be central to our analysis.
{\em For the rest of this paper, we use ``formula" to mean ``irredundant formula"}
(but we will still sometimes retain the ``irredundant" for emphasis).

\medskip
We feel sure that the analogues of Theorems \ref{SATthm} and \ref{Irr}
hold for
any fixed $k$ in place of 3; that is (with
$I_k(n)$ the number of
irredundant $k$-SAT formulas of $n$ variables),
\glossary{name={$I_k(n)$},description={The number of irredundant $k$-SAT formulae of $n$ variables.},sort=I}
we should have
\begin{conj}
\label{SATconj}
For each fixed $k$,
$G_k(n)\sim I_k(n)
\sim 2^{n+\C{n}{k}}.$
\end{conj}

\nin
While we do think it should be possible to prove this along the
present lines, the best we can say for now is that
our argument can probably be generalized to reduce Conjecture \ref{SATconj}
to a finite problem for any given $k$;
see the remarks following Corollary \ref{step2cor}.
For example,
at this writing we are pretty sure we could do $k=4$; but as this doesn't
contribute anything very interesting beyond what's needed for $k=3$,
it seems not worth adding to the present, already very long
argument.

On the other hand, if we retreat to $k=2$ then
much of the present proof evaporates---in particular hypergraph
regularity becomes ordinary Szemer\'edi regularity---leaving perhaps the
easiest verification of \eqref{allen} to date.
(Of course---if one cares---anything based on regularity must give
far slower convergence than the argument of \cite{IK}.)

\medskip
From now on we will be concerned only with the case $k=3$,
and will say ``clause" for ``3-clause,"
``formula" for ``(irredundant) 3-SAT formula," and so on.
\glossary{name={clause},description={Shorthand for {\em 3-clause}.}}
\glossary{name={formula},description={Shorthand for {\em irredundant 3-SAT formula}.}}
Let us try to say what we can
about the proof at this point.
The argument proceeds in two phases.
The first of these---which, incidentally,
gives the asymptotics of
$\log I(n)$, though the proof doesn't need to say this---is
based on the Hypergraph Regularity Lemma (HRL) of
P. Frankl and V. R\"odl \cite{FR}, a pioneering extension
to 3-uniform hypergraphs of
the celebrated (graph) Regularity Lemma of E. Szemer\'edi
\cite{Szemeredi}.
(See e.g. \cite{Rodl-Skokan}, \cite{Gowers}
for more on the spectacular recent developments on this topic.)

A mild adaptation of some of the material in \cite{FR}
shows that each
irredundant $\cee$ is ``compatible" with some
``extended partition" $\pee^*$
(defined in Section \ref{Sketch}).
On the other hand we show---this is Lemma \ref{ML1},
the upshot of this part of the argument---that the set of $\cee$'s
compatible with $\pee^*$ is small unless $\pee^*$
is ``coherent."
Since the number of $\pee^*$'s is itself negligible
relative to what we are aiming at,
this allows us to restrict our attention to $\cee$'s
compatible with coherent $\pee^*$'s.

Coherence of $\pee^*$ turns out to imply that there is
some $z\in \{0,1\}^n$
so that for {\em any} $\cee$ compatible with $\pee^*$
{\em every} witness for $\cee$
mostly agrees (in the obvious sense) with $z$.
Once we have this we are done with $\pee^*$ and the HRL,
and, in the second phase, just need to bound
the number of $\cee$'s admitting a $z$ as above, so for example
the number of $\cee$'s for which every witness is at least 99\% zeros
(note we expect that a typical such $\cee$
uses mostly {\em positive} literals.)
While this can presumably be handled as a stand-alone statement,
we instead give a recursive bound (see \eqref{I*}) that includes
minor terms involving earlier values of $I$.

The paper is organized as follows.
Section \ref{Hypergraph} fills in what we need from
hypergraph regularity.
Once we have this
we can, in Section \ref{Sketch}, make
the preceding mumble concrete and complete the proof of
Theorem \ref{Irr} assuming various supporting results.
These are proved in the remaining sections:
after some preliminaries in Section \ref{Basics},
Sections \ref{Configurations} and \ref{Coherence}
implement the first part of the above sketch (proving
Lemma \ref{ML1}); the easy
Section \ref{Witnesses} then produces the above-mentioned
$z$ associated with a coherent $\pee^*$;
and the final part of the argument (proving \eqref{I*})
is carried out in
Section \ref{Phase2}.

\bn
{\bf Usage}

Throughout the paper we use $\log$ and $\exp$ for
$\log_2$ and $\exp_2$, and $H$ for binary entropy.
We use ``$x= 1\pm y$" for ``$x\in (1-y,(1+y))$."
\glossary{name={$x=1\pm y$},description={Shorthand for $x\in (1-y,1+y)$.},sort=PM}
With the exception of \eqref{rec} (in Section \ref{Sketch})
we always assume that $n$ is large enough to support our assertions.
Following a common abuse, we usually pretend that all large numbers
are integers, and, pushing this a little, we will occasionally
substitute, e.g., ``at most $a$" for ``at most $a+1$"
in situations where the extra 1 is clearly irrelevant.

\section{Regularity}\label{Hypergraph}

In this section we recall what we need from \cite{FR} and
slightly adapt what they do to our situation.
Our notation follows theirs as much as possible.

For a bipartite graph
$G=(A\cup B,E)$, $A'\sub A$ and $B'\sub B$,
the
{\em density} of the pair $(A',B')$
is $$d(A',B')=d_G(A',B')=|E(A',B')|/(|A'||B'|)$$
\glossary{name={$d(X,Y)$},description={the density of the pair of sets $X$, $Y$.},sort=density}
\glossary{name={density of a pair of sets}}
\glossary{name={density of a bipartite graph $G$},description={$=d(A,B)$ (where $V(G)=A\cup B$).}}
(where $E(A',B')$ is the set of edges joining $A'$ and $B'$).
In particular, the {\em density} of $G$ is $d(A,B)$.
The graph $G$ (or the pair $(A,B)$)
is $\eps$-{\em regular}
if
$|d(A',B')-d(A,B)|<\eps$
for all $A' \sub A$ and $B' \sub B$
with $|A'|> \eps |A|$
and $|B'| > \eps |B|$.
\sugg{{\bf Not used, I think:}
For a given graph $G=(V, E)$ and disjoint $A, B \sub V$,
$E(A, B)$ is the set of edges joining
$A$ and $B$
and $d(A,B)=|E(A,B)|/|A||B|$ is the {\it density} of the graph
$G_{AB}:=(A\cup B, E(A,B))$.
The pair $(A,B)$ is $\eps$-regular (for $\eps>0$)
if $|d(A',B')-d(A,B)|<\eps$
for each $A' \sub A$ and $B' \sub B$
with $|A'|, |B'| > \eps |A|$.
We call a partition $V=V_0 \cup V_1 \cup \cdots \cup V_t$
{\it balanced} if $|V_0|<t$ and $|V_1|=\cdots = |V_t|$, and
$\gd$-{\it regular} if all but at most $\gd \C{t}{2}$ of the
pairs $(V_i, V_j)$ are $\gd$-regular.}

For a set $V$ write $[V]^2$ for the collection of 2-element
subsets of $V$.
\glossary{name={$(l, t, \eps_1, \eps_2)$- partition},sort=partition}
\glossary{name={$\pee$},description={See $(l, t, \eps_1, \eps_2)$-{\em partition}.},sort=P}
An $(l, t, \eps_1, \eps_2)$-{\it partition} $\pee$ of $[V]^2$
consists of an
auxiliary partition
\beq{pee1}
V=V_0 \cup V_1 \cup \cdots \cup V_t
\enq
with
$|V_0|<t$ and $|V_1|=\cdots = |V_t|=:m$, together with a system of
edge-disjoint bipartite graphs
\begin{equation}\label{pee2}
P_{\ga}^{ij}, ~~1\leq i < j \leq t, \, 0\leq \ga \leq l_{ij} \leq l,
\end{equation}
satisfying

\begin{enumerate}
\item[(a)] $\cup_{\ga=0}^{l_{ij}} P_{\ga}^{ij}=K(V_i, V_j)
~:=
\{\{x,y\}:x\in V_i,y\in V_j\} ~~ \forall\, 1\leq i < j \leq t$, and
\item[(b)] all but at most $\eps_1 \C{t}{2} m^2$ pairs $\{v_i, v_j\}$,
$v_i\in V_i$, $v_j\in V_j$, $1\leq i < j \leq t$, are edges of
$\eps_2$-regular bipartite graphs  $P_{\ga}^{ij}$.
\end{enumerate}
\nin

\medskip
\glossary{name={equitable}}
A partition $\pee$ as above is {\it equitable}
if for all but at most $\eps_1 \C{t}{2}$ pairs $i, j$, with
$1\leq i < j \leq t$, we have
$$|P_{0}^{ij}|< \eps_1 m^2$$
and
\beq{goodpair}
|d_{P_{\ga}^{ij}}(V_i, V_j)- l^{-1}|<\eps_2
~~~\forall\,
1\leq \ga \leq l_{ij}.
\enq
Note this implies
$
(1+\eps_2l)^{-1}l<l_{ij} ~~(\leq l),
$
so in fact
\beq{lijl}
l_{ij} =l
\enq
if
$\eps_2<l^{-2}$, as will be true below.

It will be convenient to refer to $V_1\dots V_t$
(but {\em not} $V_0$) as the {\em blocks of} $\pee$
and to the $P_{\ga}^{ij}$'s with $\ga>0$ as the {\em bundles of} $\pee$.
\glossary{name={block of $\pee$},description={Any $V_i$ with $i>0$.}}
\glossary{name={bundle of $\pee$},description={Any $P^{ij}_{\ga}$ with $\ga>0$.}}

\medskip
From now on
we take
$V$ to be $X_n$, our set of Boolean variables.
For the following definitions we
fix a partition $\pee$ as above.
To simplify notation
we will often use $A,B,C$ and so on for blocks of $\pee$.
\glossary{name={triad $P$ (of $\pee$)}}
A {\it triad} of $\pee$ on a triple of (distinct)
blocks $(A,B,C)$ is
$ P=P_{ABC}=(P_{AB}, P_{BC}, P_{AC})$,
with $P_{AB}$ one of the bundles of $\pee$ joining
$A$ and $B$, and similarly for $P_{AC}$ and $P_{BC}$.\footnote{This
usage differs slightly from that in \cite{FR}, in which
triads of $\pee$ may also use $P^{ij}_0$'s; the change is convenient
for us and of course does not affect Theorem \ref{HRL} (formally
it makes the theorem a bit weaker).}
\glossary{name={subtriad}}
A {\em subtriad} of such a $P$ is then
$ Q=(Q_{AB}, Q_{BC}, Q_{AC})$ with $Q_{AB}\sub P_{AB}$ and so on.
Since we are fixing $\pee$ for the present discussion,
in what follows
we will usually drop the stipulation
``of $\pee$."

\glossary{name={triangle (of a triad $P$)}}
A {\em triangle} of a triad $P$ as above is a triangle in the graph
with edge set $P_{AB}\cup P_{AC}\cup P_{BC}$
(usually designated by its set of vertices).
\glossary{name={$T(P)$},description={The set of triangles of a triad $P$.},sort=TP}
\glossary{name={$t(P)$},description={$="|T(P)"|$.},sort=tP}
We write $T(P)$ for the set of such triangles and $t(P)$ for $|T(P)|$.
Triangles of a subtriad $Q$ and $T(Q)$, $t(Q)$ are defined similarly.

\glossary{name={pattern of a triad $P$}}
For a triad $P$ on blocks $A,B,C$, a {\em pattern} on
$P$ is $\pi:\{A,B,C\}\ra\{0,1\}$.
We interpret this as associating
a preferred literal, $\pi(x)$, with each (variable)
$x\in A\cup B\cup C$; thus, for example, for $a\in A$,
$\pi(a) $ is $a$ if $\pi(A)=1$ and $\bar{a}$ if $\pi(A)=0$.
We also write $\pi(a,b,c)$ for
the clause $\pi(a) \pi(b) \pi(c)
:=\pi(a) \wedge \pi(b) \wedge \pi(c)$
(where $a\in A, b\in B$ and $c\in C$); such a
clause is said to {\em belong} to $\pi$.

\mn
{\em Remark.}
Of course we could just define patterns directly on triples
of blocks, but the current definition will turn out to
be less troublesome.
Note that, as above, we will often give the blocks of triad
$P$ as an {\em ordered} triple, which allows us to write, e.g.,
$\pi = (1,1,0)$ without ambiguity.

\medskip
Now
fix an (irredundant) formula $\cee$, again regarded
as a set of clauses.  For a triad $P$ and pattern $\pi$ on $P$,
we set
\[
T_{\pi}=T^{\cee}_{\pi}=\{ \{x,y,z\}\in T(P) : \pi(x,y,z)\in \cee \},
\]
\glossary{name={$T_{\pi}$},sort=Tpi}

\nin
and
for the analogue for
a subtriad $Q$ of $P$ use
$T_{\pi}(Q)$.
\glossary{name={density of a pattern $\pi$},description={$d_{\pi}="|T_{\pi}"|/t(P)$.}}
Define
the {\it density of} $\pi$ to be
\beq{patdens}
d_{\pi}=d^{\cee}_{\pi}=|T_{\pi}|/t(P).
\enq
\glossary{name={density of a sequence of subtriads $\Q=(Q(s))_{s\in [r]}$}}
For a pattern $\pi$ on triad $P$,
integer $r$, and $r$-tuple
$\Q=(Q(1), \ldots, Q(r))$ of subtriads of $P$,
set
\[
d_{\pi}(\Q)=\frac{|\cup_{s=1}^{ r} T_{\pi}(Q(s))|}{|\cup_{s=1}^{r} T(Q(s))|}.
\]

\mn
\glossary{name={$(\gd, r, \pi)$-regular triad},sort=regular}
\glossary{name={$(\gd, r)$-regular triad},sort=regular}
We say $P$ is
$(\gd, r, \pi)$-{\em regular for} $\cee$ if for every $\Q$ as above
with $|\cup_{s=1}^{r} T(Q(s))|>\gd t(P)$,
we have $|d_{\pi}(\Q)-d_{\pi}|< \gd$, and
$(\gd, r)$-{\it regular for} $\cee$
if it is $(\gd, r, \pi)$-regular
for each of the eight patterns
$\pi$ on $P$ (and $(\gd, r)$-{\it irregular} otherwise).

\glossary{name={$(\gd, r)$-regular partition},sort=regular}
Finally,
$\pee$ is $(\gd, r)$-{\em regular for} $\cee$ if
\beq{RegPart}
\sum \{ t(P): P \mbox{ is a } (\gd,r) \mbox{-irregular triad of }\pee \} < \gd n^3.
\enq
Let us emphasize that in the above discussion,
the quantities subscripted by $\pi$, as well as the definitions of
regularity for triads and partitions, refer to the fixed $\cee$.

\begin{thm}
\label{HRL}
For all $\gd$, $\eps_1$ with $0 < \eps_1 \leq 2\gd^4$ and integers $t_0$
and $l_0$, and for all integer-valued functions $r=r(t,l)$ and decreasing
functions $\eps_2=\eps_2(l)$ with $\, 0 < \eps_2(l) \leq  l^{-1}$, there
are $T_0$, $L_0$ and $N_0$ such that any formula $\cee$ on $X_n$, with $n > N_0$,
admits a $(\gd, r)$-regular, equitable $(l, t, \eps_1, \eps_2)$-partition $\pee$
for some $t$ and $l$ satisfying $t_0 \leq t < T_0$ and $l_0 \leq l < L_0$.
\end{thm}

\nin
{\em Proof}.
This is given by the proof of Theorem 3.11 in \cite{FR}
(which is the same as the proof of Theorem 3.5 beginning on page 151), with
some minor modifications at the outset.
We just indicate what these are, omitting a couple definitions
that are obvious analogues of their counterparts above.
We use the initial equitable $(l_0,t_0,\eps_1,\eps_2(l))$-partition $\pee_0$
(which is defined without reference to any hypergraph)
to specify hypergraphs $\h_1\dots \h_8$, as follows.
Suppose the blocks of $\pee_0$ are $V_1\dots V_{t_0}$.
For $\pi=(\pi_1,\pi_2,\pi_3)\in \{0,1\}^3$
and $x\in V_i$, $y\in V_j$ and $z\in V_k$
with $i<j<k$,
set $\pi(x,y,z)=\pi_1(x)\pi_2(y)\pi_3(z)$
($=\pi_1(x)\wedge\pi_2(y)\wedge\pi_3(z)$),
where
$$
\psi_(x) =\left\{\begin{array}{ll}
x &\mbox{if $\psi_1=1$}\\
\bar{x} &\mbox{if $\psi_1=0$,}
\end{array}\right.
$$
and similarly for $\psi_2(y)$ and $\psi_3(z)$.
Then let $\pi^1\dots\pi^8$ be some ordering of $\{0,1\}^3$, and
for $s\in [8]$ and
$x,y,z$ as above, let $\{x,y,z\}\in \h_s$ if (and only if)
$\pi^s(x,y,z)\in \cee$.

The (only) point here is that by starting this way we guarantee
that clauses belonging to the same pattern in our eventual partition
will correspond to edges of the same $\h_s$:
Theorem 3.11 of \cite{FR} gives a partition $\pee$ as in our
Theorem \ref{HRL} in which regularity with respect to $\cee$
is replaced by regularity with respect to each of $\h_1\dots \h_8$
(which we will not define).
But for any triad $P$ of $\pee$ and pattern $\pi$ on $P$,
$(\gd,r,\pi)$-regularity for $\cee$ is the same as
$(\gd,r)$-regularity of $P$ (again, we omit the definition)
for the appropriate $\h_s$, and we are done.
(To be unconscionably picky, we should slightly adjust $\gd$,
since bounds corresponding to
\eqref{RegPart} for the $\h_s$'s will turn into a bound
$8\gd n^3$ for $\cee$.)\qed

\mn
{\em Final remark.}  In applying Theorem \ref{HRL} it will be convenient
to require that in fact
\beq{lijl'}
l_{ij} =l ~~~\forall\, i,j.
\enq
As noted in \eqref{lijl} this is automatically true for $i,j$
satisfying \eqref{goodpair}
(again, assuming $\eps_2<l^{-2}$ which will be true below);
while the assumption (equitability) that all but $\eps_1\C{t}{2}$
pairs $i,j$ do satisfy \eqref{goodpair} allows us to arbitrarily
modify the partitions of the remaining $K(V_i,V_j)$'s---we just
replace them with partitions satisfying \eqref{lijl'}---without
significantly affecting \eqref{RegPart}.
(So, to be overly precise, we get this very slightly strengthened
version of Theorem \ref{HRL} by applying the original with a slightly
smaller $\gd$.
Of course the message here is that pairs failing \eqref{goodpair}
are essentially irrelevant; indeed the only point of \eqref{lijl'}
is that it makes some things a little easier to say in Section \ref{Coherence}.)

\section{Skeleton}\label{Sketch}

In this section we give enough in the way of additional definitions to
allow us to state our main lemmas, and give the proof of
Theorem \ref{Irr} modulo the much longer proofs of these supporting
results.

We will soon need to say something concrete about our many parameters,
but defer this discussion to the end of the present section.
\glossary{name={extended partition $\pee^*$}}
\glossary{name={$\pee^*$},sort=P}
\glossary{name={$\R(\pee^*)$},sort=RP}
Given $\gd,\eps_1,t_0, l_0,r=r(t,l),\eps_2=\eps_2(l)$, and
associated
$T_0, L_0$
as in Theorem \ref{HRL},
define an {\em extended partition}
$\pee^*$ to consist of an
equitable $(l, t, \eps_1, \eps_2)$-partition $\pee$,
with $t\in [t_0,T_0]$, $l\in [l_0,L_0]$,
together with

\mn
(a) a set $\R(\pee^*)$ of triads of $\pee$
that
(i) includes no $P$ for which
some two blocks of $P$ violate (\ref{goodpair})
or some bundle of $P$ violates $\eps_2$-regularity, and (ii) satisfies
\beq{RegPart'}
\sum \{ t(P): \mbox{$P$ a triad of $\pee$  not in $\R(\pee^*)$}\} < 2\gd n^3
\enq

\nin
(we will mostly ignore triads not in $\R(\pee^*)$); and

\mn
(b)  a value
$\ddd_{\pi}=\ddd^{\pee^*}_{\pi}\in \{0,t(P)^{-1}\dots (t(P)-1)t(P)^{-1},1\}$
for each pattern $\pi$ on some $P\in \R(\pee^*)$.

\medskip
\glossary{name={triad of $\pee^*$},description={Any triad $P\in \R(\pee^*)$.}}
We will call the triads in $\R(\pee^*)$ the
{\em triads of} $\pee^*$.
\glossary{name={bundle of $\pee^*$}}
The {\em bundles of} $\pee^*$ are those $\eps_2$-regular
bundles $P_{\ga}^{ij}$ of $\pee$ for which the pair $\{i,j\}$
satisfies
(\ref{goodpair})
(so the bundles of $\pee$ that we allow in triads of $\pee^*$).
\glossary{name={triangle of $\pee^*$}}
A
{\em triangle of} $\pee^*$ is a triangle
belonging to
some triad of $\pee^*$.
\glossary{name={pattern of $\pee^*$},description={A pattern on
a triad of $\pee^*$ satisfying $
\ddd_{\pi}> 2\ddd_0$.}}
Say $\pi$ is a {\em pattern of} $\pee^*$ if it is a pattern on
some triad of $\pee^*$ {\em and}
\beq{dpid0}
\ddd_{\pi}> 2\ddd_0,
\enq
where $\ddd_0$ will be specified below.
\glossary{name={clause of $\pee^*$},description={A clause belonging to a pattern of $\pee^*$.}}
A {\em clause of} $\pee^*$ is then a clause
belonging to a pattern of
$\pee^*$; 
\glossary{name={$\Cl(\pee^*)$},description={The set of clauses of $\pee^*$.},sort=KP}
we use $\Cl(\pee^*)$ for the set of such clauses.

\glossary{name={compatible}}
\glossary{name={$\cee\sim \pee^*$},sort=CsimP}
Say a formula $\cee$ and $\pee^*$ are
{\em compatible} (written $\cee\sim \pee^*$) if
every triad $P$ of $\pee^*$ is $(\gd,r)$-regular for
$\cee$,
and has
$d^{\cee}_{\pi}=\ddd_{\pi}$ for each pattern $\pi$ on $P$.
It follows from Theorem \ref{HRL} that (for large enough $n$)
every $\cee$ is
compatible with some $\pee^*$.
(The extra ``2" on the right hand side of \eqref{RegPart'}
covers triangles involving pairs $\{i,j\}$ violating
(\ref{goodpair}).)
\glossary{name={feasible}}
We say $\pee^*$ is {\em feasible} if it is compatible with at least
one $\cee$ and in what follows
{\em always assume this to be the case.}
\glossary{name={$N^*(\pee^*)$},description={$="|\{\cee:\cee\sim \pee^*\}"|$.},sort=NP}
Set
$$N^*(\pee^*) =|\{\cee:\cee\sim \pee^*\}|.
$$
We use $N^*$ here because we will later work mostly with
\glossary{name={$\N(\pee^*)$},description={$=\{\cee\cap \Cl(\pee^*):\cee\sim \pee^*\}$.},sort=NP}
\glossary{name={$N(\pee^*)$},description={$="|\N(\pee^*)"|$.},sort=NP}
\[
\mbox{$\N(\pee^*)= \{\cee\cap \Cl(\pee^*):\cee\sim \pee^*\}~$ and
$~N(\pee^*) =|\N(\pee^*)|$.}
\]

\glossary{name={proper triad}}
\glossary{name={$\pi_P$},description={The unique pattern supported by a proper triad $P$.},sort=PiP}
Say a triad $P$ of $\pee^*$ is {\em proper} if
it supports
a unique pattern of $\pee^*$---{\em always denoted}
$\pi_{_P}$---and $\ddd_{\pi_{_P}} >1/3$.
\glossary{name={agreement of $f$ and $P$}}
Say $f:\{\mbox{blocks of $\pee$}\}\ra \{0,1\}$
and $P$
{\em agree}
if $P$ is proper and $\pi_{_P}(A)= f(A)$ for
each block $A$ of $P$.
\glossary{name={coherent extended partition}}
Finally,
say $\pee^*$ is {\em coherent} if there is an
$f$ as above such that
(with $\gz_2$ discussed below)
\beq{coherence}
\mbox{{\em all but at most $\gz_2 \C{t}{3}l^3$ triads of $\pee^*$
agree with $f$}}.
\enq

The longest part of our argument is devoted to proving,
for $\ctwo$ and all of the preceding parameters
as described below,
\begin{lemma}\label{ML1}
If
\beq{bigP}
\mbox{$\log N^*(\pee^*) > (1-\ctwo)\C{n}{3}$}
\enq
then $\pee^*$ is coherent.
\end{lemma}

The argument then proceeds as follows.
Fix $\gd$, $\eps_1$, $t_0$, $l_0$, $r$, $\eps_2$
(again, see below for settings;
note $r$ and $\eps_2$ are functions).
As noted above,
Theorem {\ref{HRL}} implies that each (irredundant) $\cee$
is compatible with some extended partition $\pee^*$.
The number of possiblities for $\pee^*$ is,
for large enough $n$, less than (say)
$\exp[(\log L_0)n^2]$.
(There are, very crudely, at most: $T_0^n$ choices
for the partition $\{V_i\}$;
$\exp[(\log L_0)\C{n}{2}]$ for
the bundles $P^{ij}_{\ga}$;
and $\exp[(1+8\log m^3) \C{T_0}{3}L_0^3]$
for $\R(\pee^*)$ and the $\ddd_{\pi}$'s.)
Combining this with Lemma \ref{ML1} we have,
for any constant $c'<\ctwo$ and large enough $n$,
\begin{cor}\label{compatcor}
All but at most $\exp[ (1-c')\C{n}{3}]$
irredundant $\cee$'s satisfy
\beq{CsimP}
\mbox{$\cee\sim \pee^*$ for some {\em coherent} $\pee^*$.}
\enq
\end{cor}
We next need a bound on the number of $\cee$'s
that do satisfy (\ref{CsimP}).
\glossary{name={multiplicity}}
\glossary{name={$m(y)$},description={Multiplicity of $y$.},sort=multiplicity}
Define the {\em multiplicity}, $m(y)=m_{\cee}(y)$,
of the {\em literal} $y$ in $\cee$
to be the number of clauses of $\cee$ containing $y$.
\glossary{name={positive formula},description={$m(x)\geq m(\bar{x})$ for each variable $x$.}}
Say $\cee$ is {\em positive} if $m(x)\geq m(\bar{x})$
for each variable $x$.
The $\pee^*$'s will disappear from our argument once we
establish (with $\gz$ again TBA)

\begin{lemma}\label{Pw}
If $\cee$ is positive
and $\cee\sim\pee^*$ for some coherent $\pee^*$, then
\beq{zeta}
\mbox{any witness $\ww$ for any clause in $\cee$
has fewer than $\gz n$ $1$'s.}
\enq
\end{lemma}
\nin
The easy proof is given in Section \ref{Witnesses}.

\glossary{name={$\ii$},sort=I}
Write $\ii$ for the collection of (irredundant) positive $\cee$'s satisfying
\eqref{zeta}.
According to Corollary \ref{compatcor} and Lemma \ref{Pw}
we have
\beq{II*}
\mbox{$I(n) < \exp [(1-c')\C{n}{3}] + 2^n|\I^*|.$}
\enq
In Section \ref{Phase2} we will show, for
large enough $n$ and an appropriate positive constant $\cc$,
\begin{eqnarray}
\mbox{$|\ii|$} &<& 2^{\C{n}{3}} +
\mbox{$\exp[(1-\cc)\C{n}{2}]I(n-1)$} \nonumber\\
&&~~~\mbox{ $+ \exp[(1-\cc)3\C{n}{2}]I(n-3)
+ \exp[\C{n}{3}-\cc n] $.}~~~~\label{I*}
\end{eqnarray}

The proof of
Theorem \ref{Irr} is then completed as follows.
Combining (\ref{II*}) and (\ref{I*}) and
setting $B(n) =2^{n+\C{n}{3}}$, we have
(again, for large enough $n$)
\begin{eqnarray}
I(n)&<&
\mbox{$(1+\exp[-\cc'n])B(n)
+\exp[(1-\cc')\C{n}{2}]I(n-1)$}
\nonumber\\
&&~~~~~~~~~~~~~~~~~~~\mbox{ $
+ \exp[(1-\cc')3\C{n}{2}]I(n-3) $}\label{I(n)}
\end{eqnarray}
(where the change from $\cc$ to $\cc'$ takes care of some factors
$2^n$ and allows us to absorb the first term on the r.h.s. of (\ref{II*})
in the term $\exp[-\cc'n]B(n)$).

We show by induction that (\ref{I(n)}) implies
that, for some
constant $\gD$ and {\em all} $n$,
\beq{rec}
I(n)\leq (1+\Delta\cdot 2^{-\cc'n})B(n)
\enq
(which proves Theorem \ref{Irr}).

For (\ref{rec}),
choose $n_0$ large enough so that (\ref{I(n)})
holds for $n\geq n_0$, and then
choose $\Delta >2$ (say)
so that (\ref{rec}) holds for $n\leq n_0$.
Assuming
(\ref{rec}) holds up to $n-1$ ($\geq n_0$), we have
(omitting the little calculation for the second inequality)
\begin{eqnarray}
I(n)- B(n)&<& \mbox{$2^{-\cc'n}B(n)
+\exp[(1-\cc')\C{n}{2}](1+\gD2^{-\cc'(n-1)})B(n-1)$}
\nonumber\\
&&~~~~~~~~~~~~~~~~\mbox{$
+ \exp[(1-\cc')3\C{n}{2}](1+\gD2^{-\cc'(n-3)})B(n-3) $}
\nonumber\\
&<&
\mbox{$\{2^{-\cc'n}
+\exp[(-\cc'\C{n}{2}+n](1+\gD2^{-\cc'(n-1)})\}B(n) $}
\nonumber
\end{eqnarray}
This gives (\ref{rec}) for $n$.
\qed

\mn
{\bf Parameters}

Before proceeding we should say something about relations
between parameters.
Our task in Section \ref{Phase2} is to prove (\ref{I*})
with {\em some} positive $\cc$.  This requires
an upper bound on the $\gz$ produced by Lemma \ref{Pw}
(see \eqref{cbds} and \eqref{xibd}, which involve some
additional parameters),
which in turn, {\em via} Lemma \ref{Pw},
forces $\gz_2$ in \eqref{coherence} to be small
(namely it should satisfy \eqref{param1}).

Of course for Lemma \ref{ML1} to hold,
we then need $\ctwo$ to be small.
Specific requirements (which, for whatever it's worth,
can be satisfied e.g. with $\ctwo $ some smallish multiple of
$\gz_2^6$)
are given in Section \ref{Coherence} (see \eqref{gz1gz2}-\eqref{gz1c2c1}).
These again involve some auxiliaries, mainly $\gz_1$ and $\cone$,
which play roles in Lemma \ref{step3} analogous to those
of $\gz_2$ and $\ctwo$ in Lemma \ref{ML1}.
(The subscripts are arranged in this way because we think of
$\gz_1$ and $\cone$ as appearing earlier in the argument,
Lemma \ref{step3} being the final intermediate step in the proof
Lemma \ref{ML1}.)

We then take $\ddd_0$ to
be small compared to $\ctwo$ (the smallest of the preceding parameters),
and all of
$\gd,\eps_1,t_0^{-1},l_0^{-1}$ small compared to $\ddd_0$
(where ``small" means small enough to support our arguments;
here we won't spell out the requirements, but it will be clear
as we proceed that there is no difficulty in arranging this).
Though unnecessary,
it will be slightly convenient to set
\beq{gdt}
\gd=t_0^{-1}=l^{-1}
\enq
(but we retain the names to preserve the flavor of Theorem \ref{HRL}).
Finally, we take
$r $ ($=r(t,l)$) $=l^6$ and $\eps_2$ ($=\eps_2(l)$) $=l^{-40}$.
(The value of $r$ is needed in Section \ref{Configurations}
and then the rather severe value of $\eps_2$ is dictated by
Lemma \ref{Lseq'} (whose $h$ will eventually turn into $r$).)

We will use the usual asymptotic notation
$\ga=O(\gb)$, even when $\ga$ and $\gb$ are themselves
(usually very small) constants, the interpretation being
that $\ga < C\gb$ for some $C$ that
could be fixed in advance of any of our arguments.
But we will also sometimes use inequalities with explicit constants,
where this seems to make the exposition clearer.

\section{Basics}\label{Basics}

Here we collect some general observations,
first (Section \ref{Decency}) for regular {\em graphic}
partitions, and then (Sections \ref{Triads} and \ref{MB})
for feasible $\pee^*$'s.
These will be used in establishing, in Section \ref{Configurations},
limits on legal
configurations of patterns,
the technical basis for the proof of
Lemma \ref{ML1}.
We begin with some

\mn
{\em Conventions.}

From this point through the end of Section \ref{Coherence} we fix
a feasible $\pee^*$ (for which we will eventually prove
Lemma \ref{ML1}) together with some $\cee\sim \pee^*$.
Triads, clauses and patterns are then understood to be
triads, clauses and patterns {\em of} $\pee^*$, and we will drop
the latter specification.

As noted above, Section \ref{Decency} deals only with
graphic aspects of $\pee^*$, so does not really require feasibility.
Most of the remaining sections do require
feasibility, and it is to make use of this assumption
that we need $\cee$; that is,
we are not really interested
in $\cee$ itself at this point, but only in the implications
for $\pee^*$ that can be derived from its compatibility with $\cee$.
For the duration of this discussion
(that is, through Section \ref{Coherence}),
notation involving patterns (e.g. $T_{\pi}$)
and choices of witnesses will always
refer to $\cee$.

We will also assume,
here and in Section \ref{Configurations}, that
we have fixed a bundle $P_{\ga}^{ij}$ {\em of} $\pee^*$ for any
pair of blocks $\{V_i,V_j\}$ used by some triad involved in our
discussion; thus if two of these triads share a pair of blocks, then they
use the same bundle from this pair.
The bundles and triads under discussion may then by specified
by their blocks:
for simplicity we will usually rename blocks $A,B,C,\ldots$
and use $P_{AB}$ for the (fixed) bundle joining $A$ and $B$
and $P_{ABC}$ for the triad on $\{A,B,C\}$.
To avoid repeated specification, {\em we will always take $a, a_i\in A$
and so on.}

We will also adopt the following abusive but convenient
notation.
For blocks $A,B,C$ and $X\sub A$, $Y\sub B$, $Z\sub C$, we will write
$XY$ for the set of edges of $P_{AB}$ joining $X$ and $Y$, and
$XYZ$
for the set of triangles of the subtriad $(XY,XZ,YZ)$
of $P_{ABC}$.

\glossary{name={$Y(x_1 \dots x_k)$},sort=I}
Finally,
for a graph $G$ on $V$,
$Y\sub V$ and $x_1\dots x_k\in V\sm Y$, we set
$
Y(x_1\dots x_k)= \{y\in Y: y\sim x_i ~\forall\, i\in [k]\}
$
(where, as usual, $x\sim y$ means $xy\in E(G)$).

\mn
\subsection{Decency}\label{Decency}

We first need a few easy consequences of {\em graphic} regularity,
beginning with the following basic (and standard) observation
(see e.g. Fact 1.3 in \cite{KomSim}).

\begin{prop}\label{Pbasic}
If $(A,B)$ is $\eps$-regular with density $d$,
then for any $B'\sub B$ of size at least $\eps |B|$,
$$
|\{a\in A: |B'(a)|\neq (d\pm \eps)|B'|\}| <2\eps |A|.
$$
\end{prop}

\medskip
Now suppose that $Y_1\dots Y_k$ are (distinct) blocks of $\pee^*$
and, for $1\leq i<j\leq k$, $P_{ij}$ is a bundle of $\pee^*$
joining $Y_i$ and $Y_j$
(so in particular $P_{ij}$ is
$\eps_2$-regular with density $l^{-1}\pm \eps_2$).
\glossary{name={decent}}
For distinct $x_1\dots x_s\in \cup Y_i$
and $Y_j(x_i:i\in I)$ defined by the $P_{ij}$'s,
say $\{x_1\dots x_s\}$ is {\em decent}
(with respect to $Y_1\dots Y_k$ and the $P_{ij}$'s, but
we will drop this specification when the meaning is clear)
if for all $I\sub [s]$,
$$
|Y_j(x_i:i\in I)|= (l^{-1}\pm 2\eps_2)^{|I|}m ~~~ ( = (1\pm \es)^sml^{-s})
$$
whenever the left side is defined; that is, whenever
$x_i\not\in Y_j ~\forall\, i\in I$.

\mn
The next easy observation is
similar to, e.g.,
\cite[Fact 1.4]{KomSim}.

\begin{prop}\label{Ldecency}
With notation as above, if s is fixed and $\{x_1\dots x_s\}$
is decent,
then
for any $u\in [k]$,
$$|\{x\in Y_u:\mbox{$\{x_1\dots x_s,x\}$ is indecent}\}| < 2^{s+1}k\eps_2 m.$$
\end{prop}

\nin
(Actually we will always have $k\leq 4$, but it is no harder
to give the general statement.  In fact
$s$ need not be fixed:  we just need $(l^{-1}- 2\eps_2)^s>\eps_2$.
It may also be worth noting that the constant
$2^{s+1}k$ can always be improved; but all we ever really
need from Proposition \ref{Ldecency} is
a bound of the form $O(\eps_2m)$, so there's no
reason to be careful here.)

\mn
{\em Proof.}
If $x\in Y_u$ and $\{x_1\dots x_s,x\}$ is indecent, then there
are $j\in [k]\sm \{u\}$ and $I\sub [s]$
such that
$x_i\not\in Y_j$ $\forall\, i\in I$ and
$$
|Y_j(x)\cap Y_j(x_i:i\in I)|\neq (l^{-1}\pm 2\eps_2)|Y_j(x_i:i\in I)|.
$$
But by Proposition \ref{Pbasic}
(using $|Y_j(x_i:i\in I)|>(l^{-1}- 2\eps_2)^sm >\eps_2m$),
the number of such $x$'s for a given $j$ and $I$ is less than
$2\eps_2 m$.\qed

In line with the conventions given at the beginning
of this section, we will in what follows always assume that ``decency"
refers to the set of blocks under discussion,
and will tend to
drop the specification
``with respect to $Y_1\dots Y_k$."

\medskip
From now until the end of Section \ref{Triads} we work with
blocks $A,B,C$, employing the conventions discussed earlier
and setting $P=P_{ABC}$.
The following definitions
are given with $A,B,C$ in particular roles,
but of course are meant to also apply when these roles are permuted.
\glossary{name={$L(a)$},description={$=L_P(a)=\{bc: \, \{a,b,c\}\in T(P)\}$.},sort=La}
\glossary{name={$L(ab)$},description={$=L_P(ab)=\{c: \, \{a,b,c\}\in T(P)\}$.},sort=Lab}
Set (for $a\in A$)
$$
L(a)=L_P(a)=\{bc: \, \{a,b,c\}\in T(P)\}
$$
($L$ for ``link"),
and, similarly, for an edge $ab$,
$$
L(ab)=L_P(ab)=\{c: \, \{a,b,c\}\in T(P)\}
$$
(where, recall, we assume $a\in A$ and so on).

The next proposition, in which decency is with respect to
$A,B,C$, is immediate from the definitions
\begin{prop}\label{Ptriangles}
{\rm (a)}  If a is decent then
$|L(a)|=(1\pm 2\eps_2l)^3m^2l^{-3}$;

\mn
{\rm (b)}
If $ab$ is a decent edge, then
$|L(ab)|=(1\pm 2\eps_2l)^2ml^{-2}$;

\sugg{
\mn
{\rm (c)}
$
~(1-7\eps_2 l) m^3 l^{-3}  <t(P) <
(1+3\eps_2 l^2)m^3 l^{-3} .
$
}
\end{prop}

Finally, we need to say something about triangle counts
(compare e.g. \cite[Fact A, p. 139]{FR}):
\begin{prop}\label{Ptriangles'}
If $X,Y,Z$ are subsets of A,B,C (resp.) with each of $|Y|,|Z|$
at least $(1-\es)^{-1}\eps_2lm$,
then
$$
(1-\tfrac{2\eps_2m}{|X|})(1-\es)^3|X||Y||Z|l^{-3}<
|XYZ| < |X||Y||Z|l^{-3} + 5\eps_2m^3.
$$
In particular,
$$
(1-7\eps_2 l)m^3l^{-3}<t(P)<
(1+5\eps_2 l^3)m^3l^{-3}
$$
\end{prop}

\mn
{\em Proof}.
Lower bound:
There are at least
$|X|-2\eps_2 m = (1-\tfrac{2\eps_2m}{|X|})|X| $
$a$'s in $X$ with
$|Y(a)|>(1-2\eps_2 l)|Y|l^{-1}$ and $|Z(a)|>(1-2\eps_2 l)|Z|l^{-1}$,
and for each of these $a$'s we have
(now fully using the lower bounds on $|Y|$ and $|Z|$)
$|Y(a)Z(a)| > (1-2\eps_2 l)|Y(a)||Z(a)|l^{-1}$.

Upper bound:
There are at most $2\eps_2 m$ $a$'s with $|Y(a)|>(1+2\eps_2 l)|Y|l^{-1}$
or $|Z(a)|>(1+2\eps_2 l)|Z|l^{-1}$ (or both), while for {\em any} $a$ we have
$|Y(a)Z(a)| < \max\{(1+\es)|Y(a)||Z(a)|l^{-1}, \eps_2m^2\}$.
This gives (crudely)
$$|XYZ| < (1+\es)^3|X||Y||Z|l^{-3}+4\eps_2 m^3.$$\qed

\subsection{Triads}\label{Triads}

We continue to work with blocks $A,B,C$ and $P=P_{ABC}$,
and now fix a pattern $\pi$ on $P$.
Note in particular that ``decency" in this section is
with respect to these three blocks (and $P$).
\glossary{name={$L^{\pi}(a)$},description={$=\{bc: \, \{a,b,c\}\in T_{\pi}\}$.},sort=Lpia}
\glossary{name={$L^{\pi}(ab)$},description={$=\{c: \, \{a,b,c\}\in T_{\pi}\}$.},sort=Lpiab}
Set (e.g.)
\[
L^{\pi}(a)=\{bc: \, \{a,b,c\}\in T_{\pi}\}
\]
(where, recall, $T_{\pi}$ is $T_{\pi}^{\cee}$ for our fixed $\cee$)
and, for an edge $ab$,
\[
L^{\pi}(ab)=\{c: \, \{a,b,c\}\in T_{\pi}\}.
\]

\glossary{name={good vertex}}
Say $a$  is {\it good} for $\pi$---or for now simply {\em good}---if,
with $\gd_1=\sqrt{\gd}$,
\begin{enumerate}
\item [(i)] $a$ is decent, and
\item [(ii)]
for any $B_1\dots B_r\sub B(a)$ and $C_1\dots C_r\sub C(a)$,
if the edge sets  $G_s:= B_s C_s$
satisfy
$|\cup_{s=1}^{r} G_s|>\gd_1 m^2 l^{-3}$, then
\beq{gooda}
|L^{\pi}(a) \cap (\cup_{s=1}^{r} G_s)| = (d_{\pi}\pm \gd)|\cup_{s=1}^{r} G_s|.
\enq
\end{enumerate}
(Note that (ii) implies the formally more general statement
where the number of $B_i$'s and $C_i$'s is
{\em at most} $r$, since we can
add some empty sets to the list.)

\glossary{name={nice vertex}}
For a good $a$, say $b\in B(a)$ is {\em nice to} $a$
(with respect to $\pi$, but again we'll drop this specification) if
$\{a,b\}$ is decent
and
\beq{deg_pattern_ab}
|L^\pi(ab)| = (d_{\pi}\pm 2\gd)m l^{-2}
\enq
\glossary{name={good edge}}
An edge $ab$ is then {\it good}
if $a$ and $ b$ are
good and nice to each other.
\glossary{name={good triangle}}
\glossary{name={great vertex/edge/triangle}}
A triangle $\{a,b,c\}$ is {\it good} if its
edges are all good
and {\it great} if it is good and belongs to $T_{\pi}$.
Finally,
we say a vertex is {\it great}
if it belongs to at least
$\ddd_0 m^2l^{-3}$ great triangles and an edge is
{\it great} if it belongs to at least $\ddd_0 m l^{-2}$
great triangles.

\medskip
Let $\gd_2=4\eps_2l+3 \gd_1$, $\gd_3=12\eps_2 +3\gd_1  $,
$\gd_4=114\eps_2l^3 + 4\gd_2+4\gd_3$, and
$\gc = 2\gd_4/\ddd_0$ (= $\Theta(\sqrt{\gd}/\ddd_0)$).
We will use these ugly expressions in the statement and proof
of the next lemma, but will then immediately pass to the
relaxed version, Corollary \ref{goodnesscor},
at which point $\gd_2,\gd_3,\gd_4$ will disappear from the discussion.

\begin{lemma}\label{goodness}
{\rm (a)}
At least $(1-\gd_2) m$ vertices of $A$ are good.

\mn
{\rm (b)} If $a $ is good, then
$|\{b\in B(a):\mbox{$b$ is not nice to $a$}\}|<\gd_3ml^{-1}$;
thus at least $(1 -2\eps_2l-\gd_3) m l^{-1}$ vertices of $ B(a)$ are nice to $a$.

\mn
{\rm (c)}
At most $\gd_4 m^3 l^{-3}$ members of $T(P)$ are not good.
It follows that $T(P)$ contains at least
$(1-7\eps_2 -\gd_4)m^3 l^{-3}$ good triangles and at least
$(d_{\pi}-7\eps_2l -\gd_4)m^3l^{-3}$ great triangles.

\mn
{\rm (d)} At least $(1-\gc)m$ vertices of $A$ are great, and
at least $(1-\gc)(m^2 l^{-1})$ edges of $P_{AB}$ are great.
\end{lemma}

\begin{cor}\label{goodnesscor}
{\rm (a)}
At least $(1-\gc) m$ vertices of $A$ are good.

\mn
{\rm (b)} If $a $ is good, then
$|\{b\in B(a):\mbox{$b$ is not nice to $a$}\}|<\gc ml^{-1}$, and
at least $(1 -\gc) m l^{-1}$ vertices of $ B(a)$ are nice to $a$.

\mn
{\rm (c)}
At most $\gc m^3 l^{-3}$ members of $T(P)$ are not good.
At least
$(1-\gc)m^3 l^{-3}$ triangles of $T(P)$ are good,
and at least
$(d_{\pi}-\gc)m^3l^{-3}$ are great.

\mn
{\rm (d)} (Repeating:)
At least $(1-\gc)m$ vertices of $A$ are great and
at least $(1-\gc)m^2 l^{-1}$ edges of $P_{AB}$ are great.
\end{cor}

\mn
{\em Proof of Lemma} \ref{goodness}.
We use ``bad" for ``not good" and
for the proofs of (a) and (b) set $G=P_{BC}$.

\mn
(a)
By Proposition \ref{Pbasic},
at most $4\eps_2 m$ vertices of $A$ are indecent;
so failure of (a) implies that there is a set $A_0$ of at least
$ (3/2) \gd_1 m$ decent vertices
of $A$ satisfying either (i) for each $a\in A_0$
there are $B_1(a)\dots B_r(a)\sub B(a)$ and $C_1(a)\dots C_r(a)\sub C(a)$
such that, with $G_s(a) = B_s(a)C_s(a)$, we have
$|\cup_{s=1}^{r} G_s(a)|>\gd_1 m^2 l^{-3}$, and
$$
|L^{\pi} (a) \cap (\cup_{s=1}^{r} G_s(a))| <
(d_{\pi}-\gd)|\cup_{s=1}^{r} G_s(a)|,
$$
or (ii) the corresponding statement with ``$<(d_{\pi}-\gd)$"
replaced by ``$>(d_{\pi}+\gd)$."
Assuming the first (the argument for the second is identical) and setting
$G_s=\cup_{a\in A_0}G_s(a)$,
$H_s=\cup_{a\in A_0}\{ab:b\in B_s(a)\}$ and
$K_s=\cup_{a\in A_0}\{ac:c\in C_s(a)\}$,
we find that for the subtriads
$Q_s=(G_s,H_s,K_s)$ of $P$ we have
$$
|\cup_{s=1}^r T(Q_s)|=\sum_{a\in A_0} |\cup_{s=1}^{r} G_s(a)|>
|A_0| \gd_1 m^2 l^{-3} \geq \gd t(P)
$$
(using the upper bound on $t(P)$ in Proposition \ref{Ptriangles'}), while
\begin{eqnarray*}
|\cup_{s=1}^{ r} T_{\pi}(Q_s)|
&=&\sum_{a\in A_0} |L^{\pi}(a) \cap (\cup_{s=1}^{r} G_s(a))|\\
&<&
\sum_{a\in A_0} (d_{\pi}-\gd)|\cup_{s=1}^{r} G_s(a)|
=(d_{\pi}-\gd) |\cup_{s=1}^r T(Q_s)|,
\end{eqnarray*}
\nin
contradicting the $(\gd, r, \pi)$-regularity of $P$.

\bn
(b)
Since $a$ is decent, each of $|B(a)|,|C(a)|$ is at least
$(1-2\eps_2l)ml^{-1}$; in particular the second assertion
in (b) follows from the first.
By Proposition \ref{Ldecency},
$|\{b\in B(a):ab ~\mbox{is indecent}\}|<12\eps_2 m$;
so we will be done if we show that at most $3\gd_1 ml^{-1}$
$b$'s violate (\ref{deg_pattern_ab}).
Suppose instead
(e.g., the other case again being similar)
that there is $B_0\sub B(a)$ of size
at least $(3/2)\gd_1ml^{-1}$ with
$$
|L^\pi (ab)| < (d_{\pi}-2\gd) ml^{-2}   ~~~\forall\, b\in B_0.
$$
Then with
$G_1=B_0C(a)$ we have
\[
|G_1| > |B_0||C(a)|(l^{-1}-2\eps_2) > |B_0|(1-2\eps_2l)^2ml^{-2}
> \gd_1 m^2l^{-3},
\]
while
\[
|L^{\pi}(a) \cap G_1| = \sum_{b\in B_0} |L^{\pi}(ab)|
<  |B_0|(d_{\pi}-2\gd) m l^{-2} < (d_{\pi}-\gd)|G_1|,
\]
contradicting the assumption that $a$ is good.

\mn
(c)
Of the triangles $\{a,b,c\}$ of $T(P)$
at most $114\eps_2m^3$
are indecent
(by Proposition \ref{Ldecency}; the constant
is of course a bit excessive);
at most $3\gd_2(1+2\eps_2l)^3m^3l^{-3}<4\gd_2m^3l^{-3}$
are bad because  at least one of $a,b,c$
is decent but bad
(by (a) and Proposition \ref{Ptriangles}(a));
and
at most $3m(\gd_3ml^{-1})(1+2\eps_2l)^2ml^{-2} < 4\gd_3m^3l^{-3}$ are
decent but bad because one of $a,b,c$
fails to be nice to another
(by (b) and Proposition \ref{Ptriangles}(b)).
This gives the first assertion;
the second and third
then follow from Proposition \ref{Ptriangles'},
the latter since
the number of great triangles of $P$ is at least
\[
|T_{\pi}| -\gd_4 m^3l^{-3} = d_{\pi}t(P) -\gd_4m^3l^{-3} >
[d_{\pi}-7\eps_2l -\gd_4]m^3l^{-3}.
\]

\mn
(d)
Set $\eta= 7\eps_2l+\gd_4$; thus (c) says that the number of great triangles
is at least $(d_{\pi}-\eta)m^3l^{-3}$.

We first consider great vertices $a$.
A good $a$ belongs to at most
$(d_{\pi}+2\gd)m^2l^{-3}$ great triangles
(namely,
$|L^{\pi}(a)|<
(d_{\pi}+\gd)|(B(a)C(a))| <(d_{\pi}+\gd)
(l^{-1}+2\eps_2)^3m^2<(d_{\pi}+2\gd)m^2l^{-3}$).
Thus, with $s$ the number of non-great $a$'s
(note a bad vertex is in no great triangles), the
number of great triangles is at most
$$
(m-s) (d_{\pi}+2\gd)m^2l^{-3} + s\ddd_0m^2l^{-3},
$$
and combining these bounds gives (using \eqref{dpid0})
$s < (2\gd +\eta)/(d_{\pi}+2\gd-\ddd_0) m <\gc m$.

The argument for edges is similar.
A good edge belongs to at most $(d_{\pi}+2\gd)ml^{-2}$
great triangles, so if $s$ is the number of non-great $ab$'s
then the number of great triangles is at most
$
((1+\eps_2l)m^2l^{-1}-s)(d_{\pi}+2\gd) +s \ddd_0ml^{-2}$.
Again combining with (c) bounds
$s$ by roughly
$(\gd_1/\ddd_0) m^2l^{-1}$, and the (second)
statement in (c) follows since $|AB| >(1-\eps_2l)m^2l^{-1}$.
\qed

\subsection{More basics}\label{MB}

We continue to work with  $P=P_{ABC}$,
and a fixed $\pi$ on $P$.
For the next lemma we add a fourth block, say $D$,
which only appears incognito:  ``decency" in
Lemma \ref{Lseq'} means with respect to $A,B,C,D$.
\begin{lemma}\label{Lseq'}
For $\T\sub T(P)$ with $|T(P)\sm\T|<5\gc m^3l^{-3}$
and
$h$ such that $h^6\eps_2l^2<<d_{\pi}$,
there are distinct $a_i,b_{ij}$ and $c_{ij}$,
$i,j\in [h]$ satisfying

\mn
{\rm (i)}  $~\{a_i,b_{ij},c_{ij}\}\in \T$ is great for all $i,j$, and
\mn
{\rm (ii)}  any set of four of the vertices
$a_i,b_{ij},c_{ij}$ is decent.

\end{lemma}
\nin
In practice $\T$ will consist of all members of $T(P)$ avoiding
some set of pathologies that are known to be rare by the results of
Section \ref{Triads}.

\mn
{\em Proof.}
We first observe that, with $\ga = d_{\pi}$ and $\T^*$ the set
of great triples from $\T$,
we have (using Proposition \ref{Ptriangles'} and Corollary \ref{goodnesscor}(c))
\beq{TTpi}
|\T^*| > \ga t(P) -5\gc m^3l^{-3}-\gc m^3l^{-3}
> (\ga-6\gc- 7\eps_2l^3) m^3l^{-3}.
\enq

Say an edge $ab$ is {\em fine} if
$|\{c:abc\in \T^*\}|> \frac{1}{2}\ga ml^{-2}$, and
$a$ is {\em fine}
if $ab$ is fine for at least $\frac{1}{2}ml^{-1}$ $b$'s.
We assert that
\beq{fine1}
\mbox{at most $40(\gc/\ga) m$ $a$'s are not fine.}
\enq

\nin
{\em Proof of} (\ref{fine1}).
Writing $s$ for the number
of non-fine $ab$'s we find (with explanations to follow)
that $|\T^*|$ is at most
\beq{fine2a}
3\gc m^3l^{-3} +
((1+2\eps_2l)m^2l^{-1}-s)(\ga+2\gd)ml^{-2} + (1/2)s\ga ml^{-2}.
\enq
Here the first term
covers triangles on edges $ab$ that are either indecent
or for which
$|L^{\pi}(ab)|> (\ga +2\gd) ml^{-2}$.
(By Proposition \ref{Ldecency}
there are at most $O(\eps_2 m^2)$ $ab$'s of the first type,
a minor term since $\eps_2$ is much smaller than $\gc l^{-3}$.
On the other hand, $ab$ decent with $|L^{\pi}(ab)|> (\ga +2\gd) ml^{-2}$
implies that either
$a$ is bad, or
$a$ is good and $b$ is not nice to $a$; by
Corollary \ref{goodnesscor}(a) and (b), there are
essentially at most $2\gc m^2l^{-1}$ such $ab$'s;
decency
gives $|L^{\pi}(ab)|< (1+\es)^2ml^{-2}$.)
The expression $(1+2\eps_2l)m^2l^{-1}$ is an upper bound on the
number of decent edges $ab$, and the rest of (\ref{fine2a}) is
self-explanatory.

Combining (\ref{fine2a}) and (\ref{TTpi})
gives (say) $s < 19(\gc/\ga)m^2l^{-1}$.
It follows (using
$|P_{AB}|>(1-\es)m^2l^{-1}$) that for
the number, say
$u$, of fine $ab$'s, we have
\beq{ssmall}
u>(1-19\gc/\ga)m^2l^{-1}.
\enq
But we also have, with $v$ the number of non-fine $a$'s,
$$
u < 2\eps_2m^2 +(m-v)(1+\es)m l^{-1} +(1/2)vml^{-1}
<(m -v/2)ml^{-1} +4\eps_2m^2,
$$
and combining this with \eqref{ssmall} gives (\ref{fine1}).
\qqed

\medskip
We now turn to producing the sequences described in the lemma.
First, from the set of at least $(1-40\gc/\ga)m$
fine $a$'s,
choose (distinct)
$a_1\dots a_h$ such that
\beq{seq3}
\mbox{any 4-subset of the $a_i$'s is decent}.
\enq
This is possible because,
by Proposition \ref{Ldecency}, once we have
$a_1\dots a_i$,
(\ref{seq3}) rules out at most
$O(i^3\eps_2 m)$ choices for $a_{i+1}$.

Second, for $i=1\dots h$, do:
for $j=1\dots h$ choose (distinct) $b_{ij},c_{ij}$
with $a_ib_{ij}c_{ij}\in \T^*$
such that (ii) holds for all $a$'s, $b$'s and $c$'s
chosen to this point
(that is, any set of at most four vertices from
$\{a_1\dots a_h\}\cup
\bigcup\{\{b_{kl},c_{kl}\}:\mbox{$k<i$ or [$k=i$ and $l\leq j$]}\}$
is decent).
We can do this because (again using Proposition \ref{Ldecency})
when we come to $j$:
from an initial set of at least $(1/2)ml^{-1}$ $b$'s for
which $a_ib$ is fine,
at most $O(h^6\eps_2 m)$ are disallowed because they
introduce a violation of (ii)
or are equal to some earlier $b_{kl}$; and
similarly, given $b_{ij}$,
there are at least $(1/2)\ga ml^{-2}-O(h^6\eps_2 m)$ choices for $c_{ij}$.\qed

\medskip
In Section \ref{Configurations}
we will use sequences as in Lemma \ref{Lseq'}
to prove the impossibility of certain
combinations of patterns.
The underlying mechanism, provided by Lemma \ref{LXY},
is again similar to uses of $(\gd,r)$-regularity in
\cite{FR}.  We first need the elementary

\mn
\begin{prop}\label{Bonf}
If $S_1\dots S_h$ are sets of size at least $p$ with
$|S_i\cap S_j|<q ~~\forall\,i\neq j$, then for any $k\leq h$ we have
$$
\mbox{$|\cup S_i|\geq |\cup_{i=1}^kS_i|\geq kp-\C{k}{2}q.$}
$$
In particular, if $h\geq p/q$ then taking $k=p/q$ gives
$|\cup S_i|\geq p^2/(2q)$.
\end{prop}

\begin{lemma}\label{LXY}
{\rm (a)}
Suppose $X_i\sub A$ and $Y_i\sub B$, $i=1\dots h$ with $h > (\gl/\kappa)^2 l^{c+d-a-b}$
satisfy
\beq{XY1}
|X_i|>\gl ml^{-a}, ~~|Y_i|>\gl ml^{-b}~~~~\forall
\enq
and
\beq{XY2}
|X_i\cap X_j|< \kappa ml^{-c}, ~~ |Y_i\cap Y_j|< \kappa ml^{-d}~~~~\forall\,i\neq j,
\enq
where $\gl >\eps_2\max\{l^a,l^b\}$.
Then
\beq{XY3}
\mbox{$|\cup X_iY_i| > \tfrac{\gl^4}{3\kappa^2}m^2l^{c+d-2a-2b-1}.$}
\enq

\mn
{\rm (b)}
If $X_i\sub A$, $Y_i\sub B$ and $Z_i\sub C$,
$i=1\dots h>(\gl/\kappa)^3l^{d+e+f-a-b-c}$
satisfy
\[
|X_i|>\gl ml^{-a}, ~~|Y_i|>\gl ml^{-b}, ~~|Z_i|>\gl ml^{-c}~~~~\forall\,\forall\,i
\]
and
\[
|X_i\cap X_j|< \kappa ml^{-d}, ~~ |Y_i\cap Y_j|< \kappa ml^{-e},
~~ |Z_i\cap Z_j|< \kappa ml^{-f}~~~~\forall\, i\neq j,
\]
where (say) $\gl>40\eps_2\max\{l^a,l^b,l^c\}$ and
$\kappa > (20 \eps_2l^{d+e+f})^{1/3}$,
then 
\beq{XYZ}
|\cup X_iY_iZ_i| > \tfrac{\gl^6}{3\kappa^3} m^3l^{d+e+f-2a-2b-2c-3}.
\enq
\end{lemma}

\mn
{\em Remark.}  The assumptions on
$\gl$ and $\kappa$, as well as the precise expressions involving them
in \eqref{XY3} and \eqref{XYZ},
are best ignored.  In practice
both will be large compared to $l^{-1}$ ({\em a fortiori} to $\eps_2$),
so that the assumptions will be automatic and their roles
in the conclusions minor.
In some of our applications we could improve the constants
in these conclusions by using, e.g., different $\gl$'s in the
two bounds of \eqref{XY1}.

\mn
{\em Proof of }
(a).
We have (by (\ref{XY1}) and $\eps_2$-regularity)
\beq{XY4}
|X_iY_i|> (1-2\eps_2l)\gl^2 m^2l^{-a-b-1}  ~~~\forall\,i
\enq
and
\beq{XY5}
|X_iY_i\cap X_jY_j |= |(X_i\cap X_j)(Y_i\cap Y_j)|
<
(1+2\eps_2 l)\kappa^2m^2l^{-c-d-1}  ~~~~\forall\,i\neq j,
\enq
where the second inequality follows from $\eps_2$-regularity and
(\ref{XY2}) when each of
$|X_i\cap X_j|$, $|Y_i\cap Y_j|$ is at least $\eps_2m$,
and from
$|(X_i\cap X_j)(Y_i\cap Y_j)|\leq \eps_2m^2 $
otherwise.
Combining these and applying
Proposition \ref{Bonf}
(and sacrificing a factor like $3/2$ to take care of
the terms with $\eps_2l$'s)
gives (\ref{XY3}).

The proof of (b) is similar and we won't repeat the argument.
Here the lower bound on $|X_iY_iZ_i|$ corresponding to
(\ref{XY4}) and the upper bound on $|X_iY_iZ_i\cap X_jY_jZ_j|$
corresponding to \eqref{XY5}
are given by Proposition \ref{Ptriangles'}.\qqqed

\section{Configurations}\label{Configurations}

We continue to follow the conventions
given at the beginning of Section \ref{Basics}.

We will use (for example)

\begin{center}
\begin{tabular}{ c | c c c c}
 & A & B & C & D\\
  \hline			
  $\pi$ & $\gs_A$ & $\gs_B$ & $\gs_C $ & - \\
  $\pi'$ & $\tau_A$ & $\tau_B$ & - &$\tau_D $ \\
\end{tabular}
\end{center}

\nin
to mean that $\pi$ and $\pi'$ are
patterns
 on $P_{ABC}$ and $P_{ABD}$
respectively, with
$\pi(A)=\gs_A$ ($\in\{0,1\}$) and so on.
\glossary{name={legal configuration}}
A combination of patterns---called a
{\em configuration} and usually involving more than two patterns---is
{\em legal} if it can arise in a
feasible $\pee^*$.

Two configurations are {\em isomorphic} if they can be obtained
from each other by interchanging rows, interchanging columns,
and/or interchanging 0's and 1's within a column
(so by renaming blocks or patterns, or by interchanging the roles
of positive and negative literals within a block).
Of course legality is an isomorphism invariant.

\medskip
This long section is devoted to showing
illegality of certain configurations in a feasible
$\pee^*$.
To use the feasibility assumption we will (of course) fix
some $\cee\sim \pee^*$
and then, as usual, our notation (e.g. $L^{\pi},$ $T_{\pi}$,
witnesses)
refers to $\cee$.
We will make repeated use of Lemmas \ref{Lseq'} and  \ref{LXY},
 always with $h = r $ ($=l^6$),
$\gl= \ddd_0$, and $\kappa \approx 1$.
Usefulness of the bounds \eqref{XY3} and \eqref{XYZ}
then requires several lower bounds on $\ddd_0$, the strongest of which is
\beq{d0}
\ddd_0^8>10\gd.
\enq

Most of our configurations will involve
four blocks,
but we begin with a pair of patterns using just three,
say $A,B,C$, and abbreviate $P_{ABC}=P$.

\begin{lemma}\label{LABC}
Any two patterns for $P$ differ on at most
one of $A$,$B$,$C$.
\end{lemma}

\begin{cor}\label{CABC}
There are at most two patterns on P.
\end{cor}

\nin
{\em Proof of Lemma} \ref{LABC}.

Suppose instead that the patterns
$\pi_1$ and $\pi_2$ differ on at least two of $A$, $B$ and $C$, say
(w.l.o.g.)
$\pi_1(A)=\pi_1 (B)=\pi_1(C)=1$ and $\pi_2(B)=\pi_2(C)=0$.
There are then two cases:

\bn
\begin{tabular} {c c}
\begin{tabular}{ c | c c c }
   Case 1 & A & B & C \\
  \hline			
  $\pi_1$ & 1 & 1 & 1 \\
  $\pi_2$ & 1 & 0 & 0 \\

\end{tabular}

&

\begin{tabular}{ c | c c c }
   Case 2 & A & B & C \\
  \hline			
  $\pi_1$ & 1 & 1 & 1 \\
  $\pi_2$ & 0 & 0 & 0 \\

\end{tabular}
\end{tabular}

\bn
{\it Case} 1.
According to Lemma \ref{Lseq'} we can find $a$ ($\in A$)
and disjoint pairs $(b_i,c_i)$ ($\in B\times C$) for
$i\in [\hhh]$ satisfying:

\mn
(i)  each $\{a, b_i,c_i\}$ is great for $\pi_1$;

\mn
(ii)  $a$ is good for $\pi_2$;

\mn
(iii)  each set of three of the vertices
$a$, $b_i$, $c_i$ is decent.

\mn
To see this, let $\T$ in Lemma \ref{Lseq'} consist of those
$\{a,b,c\}\in T(P)$
for which $a$ is good for $\pi_2$.
Then Proposition \ref{Ldecency} (with $s=0$),
Corollary \ref{goodnesscor}(a) and
Proposition \ref{Ptriangles}(a) give
\[
|T(P)\sm  \T| < O(\eps_2m^3) + \gc m (1+\es)^3m^2l^{-3} <5\gc m^3l^{-3}.
\]
(Of course Lemma \ref{Lseq'} gives more than what we use here.)

\medskip
Let $\ww_i$ be a witness for
$\pi_1(a,b_i,c_i) $ ($= ab_ic_i$)
and set
$$
B_i=L^{\pi_1}(a c_i)\sm\{b_i\},~~~
C_i=L^{\pi_1}(a b_i)\sm\{c_i\}.
$$
Then for each $i$ we have
\beq{BiCi'}
|B_i|,|C_i| >  \ddd_0 ml^{-2}
\enq
(since by (i) and the definition of ``pattern of $\pee^*$"
we have
$|B_i|,|C_i| > (d_{\pi_1}-2\gd) ml^{-2}-1 > \ddd_0 ml^{-2}$) and,
by the definition of ``witness,"
$$
B_i,C_i\sub \ww_i^{-1}(0) $$
which implies that
\beq{Lpi2'''}
L^{\pi_2}(a)\cap B_iC_i =\0
\enq
(since $bc \in L^{\pi_2}(a)\cap B_iC_i$ would mean that
$\ww_i$ satisfies the clause $a\bar{b}\bar{c}\in \cee$,
contradicting the assumption that $\ww_i$ is a witness for $ab_ic_i$).
On the other hand
(iii) says
$|B_i\cap B_j|,|C_i\cap C_j|<(1+\es)^3 ml^{-3}$
($\forall\,i\neq j$),
so that, in view of (\ref{BiCi'}) and \eqref{d0},
Lemma \ref{LXY}(a) gives
$
|\cup B_iC_i| > \frac{1}{3}\ddd_0^4(1+\es)^{-6}m^2l^{-3}>
\gd_1m^2l^{-3}.
$
But then (\ref{Lpi2'''})
contradicts
(ii).

\bn
{\it Case} 2.
By Lemma \ref{Lseq'} (with $\T=T(P)$)
 we can find triples
$\{a_i,b_i,c_i\}$,
$i\in [\hhh]$, satisfying:

\mn
(i)  each $\{a_i, b_i,c_i\}$ is great for $\pi_1$;

\mn
(ii)  each set of four of the vertices
$a_i$, $b_i$, $c_i$ is decent.

\medskip
Let $\ww_i$ be a witness for
$\pi_1(a_i,b_i,c_i) $ ($= a_ib_i c_i$)
and set
$$
A_i=L^{\pi_1}(b_i c_i)\sm\{a_i\},~~
B_i=L^{\pi_1}(a_i c_i)\sm\{b_i\},~~
C_i=L^{\pi_1}(a_i b_i)\sm\{c_i\}.
$$
Then for each $i$ we have
\beq{BiCi}
|A_i|,|B_i|,|C_i| > \ddd_0 ml^{-2}
\enq
and
$$
A_i,B_i,C_i\sub \ww_i^{-1}(0), $$
the latter implying
\beq{Lpi2''}
T_{\pi_2}\cap A_iB_iC_i =\0.
\enq
On the other hand Lemma \ref{LXY}(b)
with (\ref{BiCi}) and
(ii) (which implies that each of
$|A_i\cap A_j|,|B_i\cap B_j|,|C_i\cap C_j|$ is at most $(1+\es)^4 ml^{-4}$)
gives
$$
|\cup A_iB_iC_i| > \frac{1}{3}\ddd_0^6(1+\es)^{-12} m^3l^{-3}>
\gd t(P)
$$
(where the second inequality uses \eqref{d0} and the upper bound
in Proposition \ref{Ptriangles'}),
so that (\ref{Lpi2''})
contradicts
the assumption that $\pi_2$ is a pattern.\qed

\medskip
We now turn to configurations on four blocks, say $A,B,C,D$.
At one point in the argument we will need
the next result, which is contained in Lemma 4.2 of
\cite{FR} (the ``Counting Lemma").

\begin{lemma}
\label{counting}
Let $\pi_1$, $\pi_2$, $\pi_3$ and $\pi_4$ be
patterns on $P_{ABC}$, $P_{ABD}$, $P_{ACD}$ and $P_{BCD}$, respectively.
Then for any $\cee\sim \pee^*$ there are
$a\in A$, $b\in B$, $c\in C$ and $d\in D$
so that $\pi_1(a,b,c)$, $\pi_2(a,b,d)$,  $\pi_3(a,c,d)$ and $\pi_4(b,c,d)$
are all clauses of $\cee$.
\end{lemma}

glossary{name={consistent}}
Say a configuration
is {\em consistent} if any two of its patterns agree on
their common blocks.
Our main technical result
is
Lemma \ref{illegal}, which in particular
says that, up to isomorphism,
the only {\em inconsistent}
legal configuration comprised of
patterns on three distinct triads from a given set of four blocks
is

\begin{center}
\begin{tabular}{ c | c c c c}
 Conf 0& A & B & C & D\\
  \hline			
  $\pi_1$ & $1$ & $1$ & $1$ & - \\
  $\pi_2$ & $1$ & $1$ & - & $1$ \\
  $\pi_3$ & $0$ & - & $1$ & $1$
\end{tabular}
\end{center}
(To elaborate a little: any configuration of the type described is
isomorphic to some
\begin{center}
\begin{tabular}{ c | c c c c}
& A & B & C & D\\
  \hline			
  $\pi_1$ & $1$ & $1$ & $1$ & - \\
  $\pi_2$ & $1$ & $*$ & - & $1$ \\
  $\pi_3$ & $*$ & - & $*$ & $*$
\end{tabular}
\end{center}
(where
the $*$'s are 0's and 1's); and then either the $*$'s are all 1's
(and we have coherence), or the configuration is isomorphic to
Configuration 0 above or to one of the first
eight configurations of Lemma \ref{illegal},
the only slightly nonobvious case here being the isomorphism

\bigskip
\begin{center}
\begin{tabular}{c c c}
\begin{tabular}{ c | c c c c}
 & A & B & C & D\\
  \hline			
  $\pi_1$ & $1$ & $1$ & $1$ & - \\
  $\pi_2$ & $1$ & $0$ & - & $1$ \\
  $\pi_3$ & $0$ & - & $0$ & $1$ \\
\end{tabular}

~~~& $\cong$&~~

\begin{tabular}{ c | c c c c}
   Conf 4 & A & B & C & D\\
  \hline			
  $\pi_1$ & $1$ & $1$ & $1$ & - \\
  $\pi_2$ & $1$ & $0$ & - & $1$ \\
  $\pi_3$ & $0$ & - & $1$ & $0$ \\
\end{tabular}
\end{tabular}
\end{center}

\mn
gotten by interchanging the first two rows,
the last two columns, and the 0 and 1 in the
second column.)
A convenient rephrasing of the above assertion
regarding Configuration 0 (which, again, will follow from
Lemma \ref{illegal})
is
\begin{cor}\label{illegalcor}
In a legal configuration consisting
of patterns on three different triples from a set of four blocks,
no column can contain a 0, a 1 and a blank.
\end{cor}
\begin{lemma}
\label{illegal}
The following configurations are illegal.

\bigskip
\begin{tabular}{c c }
\begin{tabular}{ c | c c c c}
   Conf 1 & A & B & C & D\\
  \hline			
  $\pi_1$ & $1$ & $1$ & $1$ & - \\
  $\pi_2$ & $1$ & $1$ & - & $1$ \\
  $\pi_3$ & $1$ & - & $0$ & $0$ \\
\end{tabular}

&

\begin{tabular}{ c | c c c c}
   Conf 2 & A & B & C & D\\
  \hline			
  $\pi_1$ & $1$ & $1$ & $1$ & - \\
  $\pi_2$ & $1$ & $1$ & - & $1$ \\
  $\pi_3$ & $0$ & - & $0$ & $0$ \\
\end{tabular}
\end{tabular}
\bigskip

\begin{tabular}{c c}
\begin{tabular}{ c | c c c c}
   Conf 3 & A & B & C & D\\
  \hline			
  $\pi_1$ & $1$ & $1$ & $1$ & - \\
  $\pi_2$ & $1$ & $0$ & - & $1$ \\
  $\pi_3$ & $1$ & - & $0$ & $0$ \\
\end{tabular}
&

\begin{tabular}{ c | c c c c}
   Conf 4 & A & B & C & D\\
  \hline			
  $\pi_1$ & $1$ & $1$ & $1$ & - \\
  $\pi_2$ & $1$ & $0$ & - & $1$ \\
  $\pi_3$ & $0$ & - & $1$ & $0$ \\
\end{tabular}
\end{tabular}
\bigskip

\begin{tabular} {c c}
\begin{tabular}{ c | c c c c}
   Conf 5 & A & B & C & D\\
  \hline			
  $\pi_1$ & $1$ & $1$ & $1$ & - \\
  $\pi_2$ & $1$ & $0$ & - & $1$ \\
  $\pi_3$ & $0$ & - & $0$ & $0$ \\
\end{tabular}

&
\begin{tabular}{ c | c c c c}
   Conf 6 & A & B & C & D\\
  \hline			
  $\pi_1$ & $1$ & $1$ & $1$ & - \\
  $\pi_2$ & $1$ & $1$ & - & $1$ \\
  $\pi_3$ & $0$ & - & $1$ & $0$ \\
\end{tabular}

\end{tabular}

\bigskip

\begin{tabular} {c c}
\begin{tabular}{ c | c c c c}
   Conf 7 & A & B & C & D\\
  \hline			
  $\pi_1$ & $1$ & $1$ & $1$ & - \\
  $\pi_2$ & $1$ & $0$ & - & $1$ \\
  $\pi_3$ & $0$ & - & $1$ & $1$ \\
\end{tabular}

&

\begin{tabular}{ c | c c c c}
   Conf 8 & A & B & C & D\\
  \hline			
  $\pi_1$ & $1$ & $1$ & $1$ & - \\
  $\pi_2$ & $1$ & $1$ & - & $1$ \\
  $\pi_3$ & $1$ & - & $1$ & $0$ \\
\end{tabular}

\end{tabular}

\bigskip

\begin{tabular} {c c}
\begin{tabular}{ c | c c c c}
   Conf 9 & A & B & C & D\\
  \hline			
  $\pi_1$ & $1$ & $1$ & $1$ & - \\
  $\pi_2$ & $1$ & $1$ & $0$ & - \\
  $\pi_3$ & $0$ & $0$ & - & $1$ \\
\end{tabular}

&

\begin{tabular}{ c | c c c c}
   Conf 10 & A & B & C & D\\
  \hline			
  $\pi_1$ & $1$ & $1$ & $1$ & - \\
  $\pi_2$ & $1$ & $1$ & $0$ & - \\
  $\pi_3$ & $1$ & $0$ & - & $1$ \\
\end{tabular}

\end{tabular}

\end{lemma}

\mn
{\em Remarks}.
The full list of forbidden configurations in Lemma \ref{illegal}
is slightly
more than what we'll eventually need, but it
seems worth recording precisely what's going on here.
Though the arguments are fairly repetitive---and we will accordingly
give less detail in the later ones---we don't see a way to
consolidate.
An outlier is Configuration 8, which is easily handled by
Lemma \ref{counting} but doesn't seem susceptible to an argument
like those for the other cases.

\mn
{\em Proof of Lemma} \ref{illegal}.
Excepting those for Configurations 7 and 8, each of the following
arguments begins with a set of variables
satisfying certain properties, with existence
again given by Lemma \ref{Lseq'}.
We only discuss this for Configurations 1 and 6
(see also Case 1 of Lemma \ref{LABC}),
arguments in the remaining cases being similar to
(usually easier than) that for Configuration 1.
Note that, without further mention,
{\em we assume in each case that the specified variables are distinct.}

\mn
{\it Configuration} 1.
Let
$a$, $b_1\dots b_r$
and $c$ satisfy:

\mn
(i)  each $\{a, b_i,c\}$ is great for $\pi_1$;

\mn
(ii)  each $\{a,b_i\}$ is good for $\pi_2$;

\mn
(iii)  $a$ is good for $\pi_3$;

\mn
(iv)  each set of three of the vertices
$a$, $b_i$, $c$ is decent

\mn
(Existence:  Take $\T$ in Lemma \ref{Lseq'} to consist
of all $\{a,b,c\}\in T(P)$
for which
$ab$ and $a$ are good for $\pi_2$ and
$\pi_3$ respectively.
Corollary \ref{goodnesscor}(a,b) bounds the number of
$a$'s that are bad for $\pi_2$ or $\pi_3$ by $2\gc m$;
the number of $b$'s that are bad for $\pi_2$ by $\gc m$;
and the number
of $\{a,b\}$'s with $a$, $b$ good for $\pi_2$ but $ab$ bad
for $\pi_2$ by $\gc m^2l^{-1}$.
Thus Propositions \ref{Ldecency} and
\ref{Ptriangles} give
$$
|T(P)\sm \T| <
O(\eps_2m^3) + \gc [m^3l^{-3} +
3 m(1+\es)^3 m^2l^{-3} + m^2l^{-1}(1+\es)m l^{-2}],
$$
which is less than $5\gc m^3l^{-3}$.)

\medskip
Let $\ww_i$ be a witness for
$\pi_1(a,b_i,c) $ ($= ab_ic$)
and set
$$
C_i=L^{\pi_1}(a b_i)\sm\{c\}~~~
\mbox{and} ~~~D_i=L^{\pi_2}(a b_i)
$$
Then
$
C_i,D_i\sub \ww_i^{-1}(0) $,
implying
\beq{Lpi3}
L^{\pi_3}(a)\cap C_iD_i =\0.
\enq
On the other hand,
$$
|C_i|,|D_i| > \ddd_0 ml^{-2}
$$
(given by (i) and (ii)) and (iv)
(which bounds each of $|C_i\cap C_j|,|D_i\cap D_j| $ by
 $(1+\es)^3 ml^{-3}$ for $i\neq j$)
imply (using
Lemma \ref{LXY}) that
$
|\cup C_iD_i| > \gd_1m^2l^{-3}.
$
But then (\ref{Lpi3})
contradicts
(iii).

\bn
{\it Configuration} 2.
Choose triples
$\{a_i,b_i,c_i\}$, $i\in [\hhh]$, satisfying:

\mn
(i)  each $\{a_i, b_i,c_i\}$ is great for $\pi_1$;

\mn
(ii)  each $\{a_i,b_i\}$ is good for $\pi_2$;

\mn
(iii)  each set of three of the vertices
$a_i$,$b_i$,$c_i$ is decent.

\medskip
Let $\ww_i$ be a witness for
$\pi_1(a_i,b_i,c_i) $
($= a_ib_ic_i$)
and set
$$
A_i=L^{\pi_1}(b_i c_i)\sm\{a_i\},~~
C_i=L^{\pi_1}(a_i b_i)\sm\{c_i\}, ~~
D_i =L^{\pi_2}(a_i b_i).
$$
Then for each $i$ we have
$$
|A_i|,|C_i|,|D_i| > \ddd_0 ml^{-2}
$$
(by (i) and (ii)) and
$$
A_i,C_i,D_i\sub \ww_i^{-1}(0). $$
The latter implies
\beq{Lpi2'}
T_{\pi_3}\cap A_iC_iD_i =\0,
\enq
while the former,
with (iii) and Lemma \ref{LXY}
(using \eqref{d0} and Proposition \ref{Ptriangles'}
as in Case 2 of Lemma \ref{LABC})
gives
$
|\cup A_iC_iD_i| > \gd m^3l^{-3},
$
and these together
contradict
the assumption that $\pi_3$ is a pattern.

\bn
{\it Configuration} 3.
Choose $a$ and
pairs
$\{b_i,c_i\}$, $i\in [\hhh]$, satisfying:

\mn
(i)  each $\{a, b_i,c_i\}$ is great for $\pi_1$;

\mn
(ii)  $a$ is good for $\pi_2$ and $\pi_3$;

\mn
(iii)  each set of three of the vertices
$a$, $b_i$, $c_i$ is decent.

\medskip
Let $\ww_i$ be a witness for
$\pi_1(a,b_i,c_i) $ ($= a b_ic_i$)
and set
$$
B_i=L^{\pi_1}(a c_i)\sm\{b_i\}, ~~~
C_i =L^{\pi_1}(a  b_i)\sm\{c_i\}
$$
and
$$
D_i^{\tau}=\ww_i^{-1}(\tau)\cap D(a), ~ ~ \tau\in \{0,1\}.
$$
Then for each $i$ we have (by (i))
$$
|B_i|,|C_i| > \ddd_0 ml^{-2}
~~\mbox{and}~~
B_i,C_i\sub \ww_i^{-1}(0). $$
W.l.o.g. there are at least $h/2$
$i$'s---say those in $I$---for which
$|D_i^1| > \frac{1}{3}ml^{-1}$, so that
Lemma \ref{LXY}
(with (iii), and just using $|D_i^1| > \ddd_0ml^{-1}$) gives
$$|\cup B_iD_i^1|~
\geq ~|\cup_{i\in I}B_iD_i^1|~ > ~\gd_1m^2l^{-3}.$$
But we also have
$$
L^{\pi_2}(a)\cap B_iD_i^1
=\0,
$$
so we contradict the assumption that $a$ is
good for $\pi_2$.

\bn
{\it Configuration} 4:
Let $\{a_{ij},c_i,d_{ij}\}$, $i,j\in [\hhh]$,
satisfy

\mn
(i)  each $\{a_{ij},c_i,d_{ij}\}$ is great for $\pi_3$;

\mn
(ii)  each $c_i$ is good for $\pi_1$;

\mn
(iii)  each set of four of the
$a_{ij}$'s, $c_i$'s and $d_{ij}$'s is decent.

\medskip
Let $\ww_{ij}$ be a witness for
$\pi_3(a_{ij},c_i,d_{ij}) $ ($= \bar{a}_{ij}c_i\bar{d}_{ij}$)
and set
$$
A_{ij}=L^{\pi_3}(c_id_{ij})\sm\{a_{ij}\},~~~
D_{ij}=L^{\pi_3}(a_{ij}c_i)\sm\{d_{ij}\}
$$
and
$$
B^{\tau}_{ij}=\ww_{ij}^{-1}(\tau)\cap B(c_i), ~ ~ \tau\in \{0,1\}.
$$
Then for all $i,j$ we have
$$
\mbox{$|A_{ij}|,|D_{ij}| > \ddd_0 ml^{-2}~~~$
and
$~~~A_{ij}, D_{ij}\sub \ww_{ij}^{-1}(1) $},$$
the latter implying in particular
that
\beq{Lpi1ci}
L^{\pi_1}(c_i)\cap A_{ij}B_{ij}^1 =\0.
\enq

Suppose first that there is an $i$ for which
$|B_{ij}^1| > \frac{1}{3}ml^{-1}$ for at least $h/2$ $j$'s,
say those in $J$.
Then combining our lower bounds on
$|A_{ij}|$ and $|B_{ij}^1|$ with (iii)
and applying Lemma \ref{LXY} gives
$$
|\cup_{j\in J}A_{ij}B_{ij}^1| > \gd_1m^2l^{-3}.
$$
But then
(\ref{Lpi1ci})
contradicts
the assumption that $c_i$ is good for $\pi_1$.

We may thus suppose (at least) that for each $i$ there is {\em some}
$j(i)$ with $|B_{i,j(i)}^0 |> \frac{1}{3}ml^{-1}$.
We then drop the remaining $j$'s and relabel
$a_i=a_{i,j(i)}$, $d_i=d_{i,j(i)}$,
$\ww_i=\ww_{i,j(i)}$,
$A_i=A_{i,j(i)}$,
$D_i=D_{i,j(i)}$ and
$B_i=B_{i,j(i)}^{0}$.

Since $A_i,D_i\sub \ww_i^{-1}(1)$ and $B_i\sub \ww_i^{-1}(0)$ we have
\beq{piABC}
T_{\pi_2}\cap (\cup A_iB_iD_i)=\0  ~~~\forall\,i.
\enq
But our lower bounds on sizes
(to repeat, these are $|A_i|, |D_i|> \ddd_0 ml^{-2}$ and $|B_i|>\frac{1}{3}ml^{-1}$)
together with (iii)
imply ({\em via} Lemma \ref{LXY}; note that here
the $|A_i\cap A_j|$'s and $|D_i\cap D_j|$'s
are all at most about $ml^{-4}$)
$$|\cup A_iB_iD_i| > \gd m^3l^{-3},$$
so that (\ref{piABC}) contradicts
the assumption that $\pi_2$ is a pattern.

\bn
{\it Configuration} 5:
Let $\{a_i,b_{ij},c_{ij}\}$, $i,j\in [\hhh]$,
satisfy

\mn
(i)  each $\{a_i,b_{ij},c_{ij}\}$ is great for $\pi_1$;

\mn
(ii)  each $a_i$ is good for $\pi_2$;

\mn
(iii) each set of four of the
$a_i$'s, $b_{ij}$'s and $c_{ij}$'s is decent.

\medskip
Let $\ww_{ij}$ be a witness for
$\pi_1(a_i,b_{ij},c_{ij}) $ ($= a_ib_{ij}c_{ij}$)
and set
$$
A_{ij}=L^{\pi_1}(b_{ij},c_{ij})\sm\{a_i\},~~
B_{ij}=L^{\pi_1}(a_i c_{ij})\sm\{b_{ij}\},~~
C_{ij}=L^{\pi_1}(a_i b_{ij})\sm\{c_{ij}\}
$$
and
$$
D^{\tau}_{ij}=\ww_{ij}^{-1}(\tau)\cap D(a_i), ~ ~ \tau\in \{0,1\}.
$$
Then
$$
|A_{ij}|,|B_{ij}|,|C_{ij}| > \ddd_0 ml^{-2}
$$
and
$$
B_{ij}, C_{ij}\sub \ww_{ij}^{-1}(0) ~~~\forall\,i,j, $$
implying in particular that
\beq{Lpi2}
L^{\pi_2}(a_i)\cap (\cup_j B_{ij}D_{ij}^1) =\0.
\enq

If there is an $i$ such that
$|D_{ij}^1| > \frac{1}{3}ml^{-1}$ for at least $h/2$ $j$'s,
then
Lemma \ref{LXY} (with (iii) and our
lower bound on
$|B_{ij}|$)
gives
$$
|\cup_j B_{ij}D_{ij}^1| > \gd_1m^2l^{-3},
$$
so that (\ref{Lpi2})
contradicts (ii).

We may thus suppose that for each $i$ there is some
$j(i)$ with $|D_{i,j(i)}^0 |> \frac{1}{3}ml^{-1}$,
and relabel
$\ww_i=\ww_{i,j(i)}$,
$A_i=A_{i,j(i)}$,
$C_i=C_{i,j(i)}$ and
$D_i=D_{i,j(i)}^{0}$.
Then $A_i,C_i,D_i\sub \ww_i^{-1}(0)$ implies
$$
T_{\pi_3}\cap (\cup A_iC_iD_i)=\0  ~~~\forall\,i,
$$
while Lemma \ref{LXY} gives
$$|\cup A_iC_iD_i| > \gd m^3l^{-3},$$
contradicting the assumption that $\pi_3$ is a pattern.

\mn
{\it Configuration} 6:
Let
$c$ and the
pairs $\{a_i,b_i\}$, $i\in [\hhh]$, satisfy

\mn
(i)  $\{a_i,b_i,c\}$ is great for $\pi_1$;

\mn
(ii)
$|L^{\pi_2}(a_i b_i)\cap D(c)| >
\ddd_0 ml^{-3}$;

\mn
(iii)  $c$ is good for $\pi_3$;

\mn
(iv)  each set of three of the vertices
$a_i$, $b_i$, $c$ is decent.

\mn
(For existence we use
Lemma \ref{Lseq'} with $\T$ consisting
of all $\{a,b,c\}\in T(P)$
for which
$|L^{\pi_2}(ab)\cap D(c)| > \ddd_0ml^{-3}$ and
$c$ is good for $\pi_3$.
(In showing $\T$ is large we restrict to $ab$'s that
are good for $\pi_2$, but this is not needed once we
have existence.)

The number of $\{a,b,c\}$'s with
$ab$ bad for $\pi_2$ or $c$ bad
for $\pi_3$ is bounded, as in the argument for Configuration 1, by
$5\gc m^3l^{-3}$.
On the other hand, if
$ab$ is good for
$\pi_2$, then $\eps_2$-regularity (of $P_{CD}$) gives
$|L^{\pi_2}(ab)\cap D(c)| > (d_{\pi_2}-2\gd)(1-\eps_2l )ml^{-3}
>\ddd_0 ml^{-3}$
for all but at most $\eps_2m $
$c$'s.)

\medskip
Let $\ww_i$ be a witness for
$\pi_1(a_i,b_i,c) $ ($= a_i b_i c$)
and set
$$
A_i=L^{\pi_1}(b_i c)\sm\{a_i\}~~~
\mbox{and} ~~~ D_i=L^{\pi_2}(a_i b_i)\cap D(c).$$
Then
$
A_i,D_i\sub \ww_i^{-1}(0) $
implies
\beq{Lpi3'}
L^{\pi_3}(c)\cap A_iD_i =\0;
\enq
but
$$
|A_i|> \ddd_0 ml^{-2}~~\mbox{and}~~ |D_i| > \ddd_0 ml^{-3}
$$
(given by (i) and (ii)) and (iv) imply
(using
Lemma \ref{LXY} and (iv);
note here $|D_i\cap D_j|$ is at most about $ml^{-5}$)
$
|\cup A_iD_i| > \gd_1m^2l^{-3},
$
so that (\ref{Lpi3'})
contradicts
(iii).

\mn
{\it Configuration} 7.
For a pattern $\pi$ on $P_{ABC}$,
say $c$ is {\em good for $\pi$ relative to $d$}
if
$L^{\pi}(c)\cap A'B'\neq \0$ whenever
$A'\sub A(c,d)$ and $B'\sub B(c,d)$ are each of size at least
$\ddd_0 ml^{-2}$;
of course ``$d$ good for $\pi'$ relative to $c$"
for a pattern $\pi'$ on $P_{ABD}$ is defined similarly.

To rule out Configuration 7 it will be enough to show that there is
{\em some} $\{a,c,d\}$ that is great for $\pi_3$
and satisfies

\mn
(i) $c$ is good for $\pi_1$ relative to $d$;

\mn
(ii) $d$ is good for $\pi_2$ relative to $c$;

\mn
(iii) $\{a,c,d\}$ is decent.

\medskip
Given such a triple, choose a witness $\ww$ for $\pi_3(a,c,d)$
and set
$$A'=L^{\pi_3}(c,d)\sm \{a\} ~~~(\sub \ww^{-1}(1))$$
and
$$B^{\tau} = \ww^{-1}(\tau)\cap B(c,d), ~~~ \tau\in \{0,1\}.$$
We then have
$|A' |>\ddd_0 ml^{-2}$
(since $cd$ is good for $\pi_3$) and, w.l.o.g.,
$|B^1|> \frac{1}{2}(1-\es)^2ml^{-2}$,
contradicting (i) (since $L^{\pi_1}(c)\cap A'B^1=\0$).

For existence of $a,c,d$ as above, we
may argue as follows.
We know from Corollary \ref{goodnesscor}(c) that at least
$\ddd_0m^3l^{-3}$ triangles $\{a,c,d\}$ are great for $\pi_3$,
so just need to show that the number that fail to satisfy
(i)-(iii) is smaller than this.
The number that violate (iii) is (by Proposition \ref{Ldecency},
as usual) $O(\eps_2m^3)$.
We will bound the number of violations of (i), and of course
the same bound applies to (ii).

By Corollary \ref{goodnesscor}(a)
at most $\gc m$ $c$'s are not good for $\pi_1$.
On the other hand, we assert that if $c$ {\em is} good for
$\pi_1$ then the size of
$D':=\{d\in D(c): \mbox{$c$ is not good for $\pi_1$ relative to $d$}\}$
is $O(h\eps_2m)$.
For suppose this is false and
choose
$d_1\dots d_h\in D'$ with all triples $\{c,d_i,d_j\}$ decent.
(For existence of the $d_i$'s just note that, as in Lemma \ref{Lseq'},
the number of $d$'s that {\em cannot} be $d_{i+1}$ is at most
$O(i\eps_2 m)$; of course this is where we use the assumption
that $D'$ is large.)
For each $i\in [\hhh]$ let $A_i\sub A(c,d_i)$ and $B_i\sub B(c,d_i)$
be sets of size at least
$\ddd_0  m l^{-2}$ with
$L^{\pi_1}(c)\cap A_iB_i=\0$;
then Lemma \ref{LXY}
(using decency to guarantee that the $|A_i\cap A_j|$'s and
$|B_i\cap B_j|$'s are small)
gives $|\cup A_iB_i|> \gd_1 m^2l^{-3}$,
so that
$L^{\pi_1}(c)\cap \cup A_iB_i=\0$ says that in fact
$c$ was {\em not} good for $\pi_1$
(so we have our assertion).
Thus the number of triangles $\{a,c,d\}$ for which $\{c,d\}$ is decent but
violates (i) is at most
$$[\gc m(1+\es)ml^{-1} +
O(h\eps_2)m^2](1+\es)^2ml^{-2}< 4\gc m^3l^{-3}.$$
\qqed

\mn
{\it Configuration} 8.
As mentioned earlier,
this one doesn't seem to follow from an
argument like those above, but is an easy consequence of
Lemma \ref{counting},
according to which there are
$a,b,c,d$
such that each of
$\pi_1(a,b,c) =abc$,
$\pi_2(a,b,d) =abd$ and
$\pi_3(a,c,d) =ac\bar{d}$ belongs to $\cee$.
But this is impossible, since
a witness $\ww$ for $a b c $ must satisfy
either
$abd$ (if $\ww(d)=1$)
or $ac\bar{d}$ (if $\ww(d)=0$).

\mn
{\it Configuration} 9.
Choose $d$ and
$\{a_i,b_i\}$,
$i\in [\hhh]$ satisfying

\mn
(i)  each $\{a_i, b_i,d\}$ is great for $\pi_3$;

\mn
(ii)  each set of four of the vertices
$a_i$, $b_i$, $d$ is decent.

\medskip
Let $\ww_i$ be a witness for
$\pi_3(a_i,b_i,d) $ ($= \bar{a}_i \bar{b}_i d$)
and set
$$
A_i=L^{\pi_3}(b_i d)\sm\{a_i\}, ~~~
B_i=L^{\pi_3}(a_i d)\sm\{b_i\}
$$
and
$$
C_i^{\tau}=\ww_i^{-1}(\tau)\cap C(a_i,b_i), ~ ~ \tau\in \{0,1\}.
$$
W.l.o.g.
$|C^1_i|> \frac{1}{3}ml^{-2}$ for at least $h/2$ $i$'s.
But then $A_i,B_i\sub \ww_i^{-1}(1)$ and
$|A_i|,|B_i|> \ddd_0 ml^{-2}$ imply
$|\cup A_iB_iC^1_i| > \gd m^3l^{-3}$, so that
$$
T_{\pi_1}\cap \cup A_iB_iC^1_i
=\0
$$
contradicts the assumption that $\pi_1$ is a pattern.

\mn
{\it Configuration} 10.
Choose $a$ and $\{b_i,d_i\}$, $i\in [\hhh]$,
satisfying

\mn
(i)  each $\{a, b_i,d_i\}$ is great for $\pi_3$;

\mn
(ii)  $a$ is good for $\pi_1$ and $\pi_2$;

\mn
(iii)  each set of four of the vertices
$a$, $b_i$, $d_i$ is decent.

\medskip
Let $\ww_i$ be a witness for
$\pi_3(a,b_i,d_i) $ ($= a\bar{b}_id_i$)
and set
$$
B_i=L^{\pi_3}(a_i d_i)\sm\{b_i\}
$$
and
$$
C_i^{\tau}=\ww_i^{-1}(\tau)\cap C(a), ~ ~ \tau\in \{0,1\}.
$$
W.l.o.g.
$|C^1_i|> \frac{1}{3}ml^{-1}$ for at least $h/2$ $i$'s.
But then $B_i\sub \ww_i^{-1}(1)$ and
$|B_i|> \ddd_0 ml^{-2}$ give
$|\cup B_iC^1_i| > \gd_1 m^2l^{-3}$ and
$$
L^{\pi_1}(a)\cap (\cup B_iC^1_i) =\0,
$$
contradicting (ii).

\qed

\section{Coherence}\label{Coherence}

Here we complete the proof of Lemma \ref{ML1}.
We continue to work with a fixed feasible $\pee^*$
(so that ``triad" and so on continue to mean ``of $\pee^*$"
unless otherwise specified).
As usual in applications of regularity,
we will eventually have to
say that we can more or less ignore some minor
effects, here those associated with clauses not belonging to patterns
of $\pee^*$;
but we delay dealing with this for as long as possible
(until we come to ``Proof of Lemma \ref{ML1}" below).

\medskip
In addition to the ``auxiliary"
parameters $\gz_1$ and $\cone$ mentioned earlier
(at the end of Section \ref{Sketch}) we use $\varphi=.05$,
chosen to satisfy
\beq{phi1} \varphi < (1-H(1/3) )/2
\enq
and
\beq{phi2}  \varphi < \min \{10-\aaa-\bbb\log 3:\aaa,\bbb\in \Nn,
\aaa+\bbb\log 3 < 10\}.
\enq
We then require
\beq{gz1gz2}
\gz_1 << \gz_2^2,
\enq
meaning $\gz_1 < \eps\gz_2^2$ for a suitable small $\eps$
which we will not specify;
\beq{c1phi}
10 \cone\varphi^{-1} < (\gz_1/6)^2;
\enq
and
\beq{gz1c2c1}
\gz_1 > 2\ctwo\cone^{-1}.
\enq
(Given $\gz_2$ we may successively choose $\gz_1$, $\cone$, $\ctwo$
small enough to achieve \eqref{gz1gz2},
\eqref{c1phi}
and \eqref{gz1c2c1} respectively.)

\medskip
\glossary{name={bundle configuration}}
\glossary{name={BC},description={Bundle configuration.}}
Define a {\em bundle configuration}
(BC) of $\pee^*$ to be any
$\gb =(\gb_{ij}:\{i,j\}\in \C{[t]}{2})\in [l]^{\C{[t]}{2}}$.
\glossary{name={$I$-bundle},sort=bundle}
\glossary{name={block of $\gb$},description={A block with index in $I$, where $\gb$ is an $I$-bundle.}}
\glossary{name={$k$-bundle},description={An $I$-bundle with $"|I"|=k$.},sort=bundle}
\glossary{name={subbundle}}
\glossary{name={bundle of $\gb$}}
Similarly, for $I\sub [t]$, an $I$-{\em bundle} is some
$\gb =(\gb_{ij}:\{i,j\}\in \C{I}{2})\in [l]^{\C{I}{2}}$.
In this case we call the blocks indexed by $I$ the
{\em blocks of} $\gb$; say $\gb$ is a $k$-{\em bundle} if $|I|=k$;
and for $J\sub I$ set $\gb[J]=(\gb_{ij}:\{i,j\}\in \C{J}{2})$---a
{\em subbundle} or $|J|$-{\em subbundle} of $\gb$.
In any case we call the
$P_{\gb_{ij}}^{i,j}$'s ($i,j$ in $[t]$, $I$ or $J$ as appropriate)
the {\em bundles of} $\gb$ (or, in the last case, $\gb[J]$).
Of course
those for which $\{i,j\}$ violates \eqref{goodpair} or $P_{\gb_{ij}}^{i,j}$
is not $\eps_2$-regular are essentially irrelevant; but they are useful
for bookkeeping purposes.

\glossary{name={pattern of $\gb$},description={A pattern of $\pee^*$ supported on bundles of $\gb$.}}
\glossary{name={clause of $\gb$},description={Clause of $\pee^*$ supported on bundles of $\gb$.}}
\glossary{name={$\Cl(\gb)$},description={Set of clauses of $\gb$.},sort=Kb}
\glossary{name={$\N(\gb)$},description={$=\{\cee\cap \Cl(\gb):\cee\sim \pee^*\}$.},sort=Nb}
\glossary{name={$N(\gb)$},description={$="|\N(\gb)"|$.},sort=Nb}
The next few definitions parallel the discussion
leading to Lemma \ref{ML1}.
The {\em patterns} and {\em clauses} of a BC or $k$-bundle
$\gb$ are those patterns and clauses
of $\pee^*$ that are supported
on (bundles of) $\gb$.  We use
$\Cl(\gb)$ for the set of clauses of $\gb$
(so the set of members of $\Cl(\pee^*)$ supported on $\gb$),
$\N(\gb) = \{\cee\cap \Cl(\gb):\cee\sim \pee^*\}$ and
$N(\gb)=|\N(\gb)|$.

\glossary{name={triad of $\gb$},description={A triad of $\pee$ supported on $\gb$.}}
In contrast we will take a {\em triad of} $\gb$ to be any
{\em triad of} $\pee$ (the partition underlying $\pee^*$)
supported on $\gb$.
But note that as soon as a triad
supports a pattern it is necessarily a triad of $\pee^*$;
in particular a {\em proper} triad of $\gb$ will be a
proper triad of $\pee^*$ supported on $\gb$.

It will now also be helpful to define
\beq{hbeta}
h(\gb) = [(1 + 5\eps_2l^3)m^3l^{-3}]^{-1}\log N(\gb)
\enq
and $h(\pee^*) = [(1 + 5\eps_2l^3)m^3l^{-3}]^{-1}\log N(\pee^*)$,
the expression in square brackets being the upper bound on $t(P)$
given by
Proposition \ref{Ptriangles'} (for any triad $P$ of $\pee^*$).
This is a convenient normalization:
for a pattern $\pi$ of $\pee^*$, say on triad $P$,
the number of possibilities for
the restriction of a $\cee\sim \pee^*$ to $\pi$
is at most
\beq{tpdp}
\Cc{t(P)}{\ddd_{\pi}t(P)} < \exp [H(\ddd_{\pi})t(P)]
\enq
(recall $H$ is binary entropy),
so that the aforementioned upper bound gives
$$
h(\gb) \leq \sum \{H(\ddd_{\pi}):\mbox{$\pi$ a pattern of $\gb$}\}.
$$
For $\gb$ a given $I$-bundle, $J\sub I$, and $A\dots Z$ the blocks
indexed by $J$, we will also write $h(A\dots Z)$ for $h(\gb[J])$.

\glossary{name={coherent $k$-bundle}}
For a fixed $k$, say a $k$-bundle $\gb$
is {\em coherent}
if there is 
$f_{\gb}:\{ \mbox{blocks of } \gb\}\ra\{0,1\}$ 
such that
each triad $P$ of $\gb$ agrees with $f_{\gb}$
(which, recall, includes the requirement that $P$ be proper).
The definition for
coherence of a BC is defined is similar to that for an extended partition;
\glossary{name={coherent BC}}
precisely:  a {\em BC} $\gb$ is {\em coherent}
if
there is some
$f=f_{\gb}:\{\mbox{blocks of $\pee^*$}\}\ra \{0,1\}$ such that
\beq{coherence'}
\mbox{{\em all but at most $\gz_1 \C{t}{3}$ triads of $\gb$
agree with $f_{\gb}$.}}
\enq

\medskip
In outline the proof of Lemma \ref{ML1} goes as follows.
First, the forbidden configuration
results of Section \ref{Configurations} are used
to prove
\begin{lemma}\label{ML2}
For a 4-bundle $\gb$,
any legal configuration consisting of one pattern on each of the
four triads of $\gb$ is consistent.
\end{lemma}
\nin
(Recall consistency was defined (in the natural way) following
the statement of Lemma \ref{counting}.)

Using this and, again, the results of Section \ref{Configurations},
we obtain what we may think of as a ``local" version of
Lemma \ref{ML1}, {\em viz.}

\begin{lemma}\label{step2}
A 5-bundle $\gb$ with
\beq{big5}
h(\gb)> 10-\varphi
\enq
is coherent.
\end{lemma}

\begin{cor}\label{step2cor}
For any 5-bundle $\gb$, $h(\gb) \leq 10$.
\end{cor}

\mn
{\em Remarks.}
Note that
the analogues of
Corollary \ref{step2cor} and Lemma \ref{step2}
for {\em 4-bundles} $\gb$
(namely that
$h(\gb)$ is at most 4 and that $h(\gb)$ close to 4 implies coherence)
are not
true; rather, $h(\gb)$ can be as large as
$3\log 3$, as shown by adding the pattern $\pi_6 = (1,1,0)$ on
$(B,C,D)$ to Configuration 11 in the proof of Lemma \ref{step2} below.
It is for this reason that we need to work with 5-bundles.

For extension of the
present results from 3 to
larger $k$, it is getting to a suitable analogue of
Lemma \ref{step2} that so far requires $k$-specific treatment,
though a general argument does not seem out of the question.
Notice for example that for $k=4$, the ``5" in Lemma \ref{step2}
will become ``7," since (compare the preceding paragraph) there
can be 6-bundles $\gb$ with $h(\gb) > 15 $ ($=\C{6}{4}$).
Here one should of course substitute \cite{Rodl-Skokan} for \cite{FR},
which does not seem to cause any difficulties.
The rest of the argument (i.e. from Lemma \ref{step2} onwards) seems
to go through
without much modification.

\medskip
Once we have Lemma \ref{step2} (and Corollary \ref{step2cor})
we are done with all that's come before, and may derive
Lemma \ref{ML1} from these last two results.
A convenient intermediate step is
\begin{lemma}\label{step3}

\mn
{\rm (a)}
For any BC $\gb$, $h(\gb)\leq\C{t}{3}$.

\mn
{\rm (b)}
Any BC $\gb$ with
\beq{bigB}
\mbox{$h(\gb)> (1-\cone)\C{t}{3}$}
\enq
is coherent.
\end{lemma}

\medskip
Before turning to proofs we need some
quick preliminaries.
We first recall {\em Shearer's Lemma}
\cite{CFGS}, which we will need here and again in
Section \ref{Phase2}.
For a set $W$,
$A\sub W$ and $\f\sub 2^W$, the {\em trace}
of $\f$ on $A$ is
${\rm Tr}(\f,A)=\{F\cap A:F\in \f\}$.
For a hypergraph $\h$ on $W$---that is, a collection
(possibly with repeats) of subsets of $W$---we use, as usual,
$d_{\h}(x)$ for the degree of $x\in W$ in $\h$;
that is, the number of members of $\h$ containing $x$.
The original statement of Shearer's lemma
(though his proof gives a more general entropy version) is
\begin{lemma}\label{Shearer}
Let $W$ be a set and $\f\sub 2^W$, and let $\h$ be a hypergraph on $W$
with
$d_{\h}(v)\geq k$ for each $v\in W$.
Then
$$\log |\f|\leq \frac{1}{k}\sum_{A\in \h}\log |{\rm Tr}(\f,A)|.$$
\end{lemma}
\nin
Applications of Lemma \ref{Shearer} in the present section will
be instances of
\begin{cor}\label{Shcor}
{\rm (a)}
Suppose $3\leq k<q$; let $I$ be a $q$-subset of $[t]$ and
$\gb$ an $I$-bundle. Then
$$
h(\gb)\leq \Cc{q-3}{k-3}^{-1}\sum\{h(\gb[J]):J\in \Cc{I}{k}\}.
$$

\nin
{\rm (b)}  $~h(\pee^*)\leq l^{-\C{t}{2}+3}\sum h(\gb)$,
where
the sum runs over BC's $\gb$ (of $\pee^*$).
\end{cor}
\nin
{\em Proof.}
For (a) apply Lemma \ref{Shearer} with
$W= \Cl(\gb)$, $\f=\N(\gb)$ and $\h=\{\Cl(\gb[J]):J\in \C{I}{k}\}$.
Then
$\Tr(\f,\Cl(\gb[J]))=\N(\gb[J])$ and $d_{\h}(C)= \C{q-3}{k-3}$
for each $C\in W$, and the statement follows.

The proof of (b) is similar and is omitted.\qqqed

\medskip
We will also make some use of the following easy (and presumably
well-known) observation, whose proof we omit.
\begin{lemma}\label{component}
Any graph $G$ with $s$ vertices and
at least $(1-\ga)\C{s}{2}$ edges (where $0\leq \ga < 1/2$)
has a component of size at least $(1-\ga)s$.
\end{lemma}
\nin

Finally, we recall that (as in \eqref{tpdp}), for any $m$ and $\ga\in [0,1/2]$,
$$
\C{m}{\ga m} < \exp [H(\ga)m].
$$

\mn
{\em Proof} of Lemma \ref{ML2}
A counterexample would be a configuration of the form\begin{center}
\begin{tabular}{ c | c c c c}
& A & B & C & D\\
  \hline			
  $\pi_1$ & $*$ & $*$ & $*$ & - \\
  $\pi_2$ & $*$ & $*$ & - & $*$ \\
  $\pi_3$ & $*$ & - & $*$ & $*$ \\
  $\pi_4$ & - & $*$ & $*$ & $*$
\end{tabular}
\end{center}
(where the $*$'s are 0's or 1's),
in which we may assume (invoking isomorphism) that each column
contains at most one 0.
Since the configuration is incoherent
there is at least one 0, say
(w.l.o.g.) $\pi_1(A)=0$.  But then
Corollary \ref{illegalcor} says that
the configuration consisting of $\pi_1,\pi_2$ and $\pi_4$
is illegal
(as is the full configuration).\qqqed

\mn
{\em Proof of Lemma} \ref{step2}.

\medskip
Suppose $A,B,C$ are blocks of $\gb$, with $P$
the corresponding triad of $\gb$.
Of course $h(A,B,C)$ is zero if there is no
pattern (of $\gb$) on $(A,B,C)$,
and at most 1 if there is exactly one such pattern.
We assert that
\beq{Nabc}
h\abc \leq\log 3
\enq
in any case
(really meaning when there are exactly two patterns
on $\abc$; see Corollary \ref{CABC}).
To see this, suppose (w.l.o.g.) $\pi=(1,1,1)$ and
$\pi'=(1,1,0)$ are patterns on $\abc$
and,
for a fixed pair $a,b$ ($a\in A,b\in B$),
consider the possibilities for the links
$L^{\pi}(ab)=L_{\cee}^{\pi}(ab)$
and $L^{\pi'}(ab)=L_{\cee}^{\pi'}(ab)$
(with $\cee\sim \pee^*$).
We cannot
have $c\in L^{\pi}(ab)\cap L^{\pi'}(ab)$
unless each of
these links consists {\em only} of $c$
(since e.g. a witness for $abc'$ ($c'\neq c$) would agree with one
of $abc$, $ab\bar{c}$).
Thus $(L^{\pi}(ab),L^{\pi'}(ab))$ is either a pair of
disjoint subsets of $C(a,b)$ ($=L_P(ab)$) or two copies of the same
singleton, whence the number of possibilities for this pair is
less than $\exp_3[|C(a,b)|]+|C(a,b)|$.
This nearly gives (\ref{Nabc}) since $\sum|C(a,b)|=t(P)$;
to keep the clean expression in \eqref{Nabc} (which of course is not
really necessary), one may use
the fact that $\cee\sim \pee^*$ requires that
$\sum_{ab}|L^{\pi} (ab)|=\ddd_{\pi}t(P)$,
but we leave this detail to the reader.
(We could also get around this by slightly shrinking the coefficient of
$\log N(\gb)$ in \eqref{hbeta}.)
\qqed

\medskip
It follows, using Lemma \ref{ML2}
and Corollary \ref{CABC}, that if $A,B,C,D$ are blocks of $\gb$,
indexed by $J$ say, with
$h(\gb[J]) > 3 +H(1/3)$ ($>2 +\log 3$), then either
$\gb[J]$ is coherent or
exactly three of its triads
support patterns, and at least two of them support two
patterns.
It's also easy to see,
using Corollary \ref{illegalcor},
that if we do have the latter
possibility, say with two patterns on each of $(A,B,D)$
and $(A,C,D)$ and at least one on $(B,C,D)$, then
up to isomorphism (the set of patterns of) $\gb[J]$ contains
the configuration

\begin{center}
\begin{tabular}{ c | c c c c}
 Conf 11& A & B & C & D\\
  \hline			
  $\pi_1$ & $1$ & $1$ & - & $1$ \\
  $\pi_2$ & $1$ & $1$ & - & $0$ \\
  $\pi_3$ & $1$ & - & $1$ & $1$ \\
  $\pi_4$ & $1$ & - & $1$ & $0$\\
  $\pi_5$ & - & $1$ & $1$ & $1$
\end{tabular}
\end{center}

We next assert that if $\gb$
is incoherent (and satisfies \eqref{big5}),
then
\beq{some4}
\mbox{some 4-subbundle $\gb'$ of $\gb$ is incoherent with $h(\gb')>3 +H(1/3)$,}
\enq
so, according to the preceding discussion, contains Configuration 11.
For the assertion, notice that incoherence of $\gb$ implies incoherence of
at least one of its 4-subbundles; so if \eqref{some4} fails,
then Corollary \ref{Shcor} (and the fact that $h(\gb')\leq 4$ for a
coherent 4-bundle $\gb'$) gives
$$h(\gb)\leq\tfrac{1}{2}[4\cdot 4+ 3 +H(1/3)]<10-\varphi,$$
contradicting \eqref{big5}.

\medskip
Assume then that $\gb$ contains Configuration 11;
let $E$ be the fifth block of $\gb$;
and let $\aaa$ be the number of triads of $\gb$
that support exactly one pattern,
and $\bbb$ the number that support exactly two.
Then
$$h(\gb) \leq \aaa+\bbb\log 3,$$
implying in particular (using \eqref{big5} and \eqref{phi2}) that
\beq{aaabbb}
\aaa+\bbb \log 3\geq 10.
\enq

Corollary \ref{illegalcor} now says:
(i) there is no pattern on $\{A,B,C\}$
(since such a pattern together with (e.g.) $\pi_1$ and $\pi_4$
would violate the corollary);
(ii) there is either no pattern on $\{A,B,E\}$ or
no pattern on either of $\{A,D,E\}$, $\{B,D,E\}$
(since if $\pi$ is a
pattern on $\{A,B,E\}$ and $\pi'$ a pattern on either of
$\{A,D,E\}$, $\{B,D,E\}$, then $\pi$ and $\pi'$
together with one of
$\pi_1$, $\pi_2$ violate the corollary);
and similarly
(iii) there is either no pattern on $\{A,C,E\}$ or
no pattern on either of $\{A,D,E\}$, $\{C,D,E\}$.

It follows that $\aaa+\bbb\leq 7$, which with \eqref{aaabbb} implies
$\bbb\geq 6$, so that there is
a set of four blocks from $\{A,B,C,D,E\}$ three
of whose triads support two patterns apiece
(since if $S_1\dots S_6$ are 3-subsets of a 5-set $S$, then
some 4-subset of $S$ contains at least three $S_i$'s).
But we have already seen, in the derivation of Configuration 11,
that any configuration consisting of five of these patterns
must be isomorphic to Configuration 11, whence it
follows easily that (up to isomorphism)
$\gb$ contains Configuration 11 together with
\begin{center}
\begin{tabular}{ c | c c c c}
& A & B & C & D\\
  \hline			
  $\pi_6$ & - & $1$ & $1$ & $0$
\end{tabular}
\end{center}
The discussion in the preceding paragraph then shows
that there is either no pattern on $\{B,C,E\}$ or
no pattern on either of $\{B,D,E\}$, $\{C,D,E\}$;
and combining this with (i)-(iii) above gives $\aaa+\bbb\leq 6$,
contradicting \eqref{aaabbb}.\qed

\mn
{\em Proof of Lemma} \ref{step3}.

\mn
(a)  This is immediate from Corollaries \ref{Shcor}(a)
(with $q=t$, $I=[t]$) and \ref{step2cor}.

\mn
(b)
We first assert that
(for $\gb$ as in \eqref{bigB})
\beq{beta5}
\mbox{all but at most $10\cone\varphi                              ^{-1}\Cc{t}{5}$ 5-bundles of
$\gb$ are coherent.}
\enq
{\em Proof.}
By Lemma \ref{step2}, the number of
incoherent 5-bundles of $\gb$ is at most
$$
s:=|\{I\in \Cc{[t]}{5}:h(\gb[I]) < 10 -\varphi\}|.
$$
Thus, again using Corollaries \ref{Shcor}(a) and \ref{step2cor},
we have
$$
h(\gb) \leq  \mbox{$\C{t-3}{2}^{-1}[(\C{t}{5}-s)10 + s (10-\varphi)]$}
=  \mbox{$\C{t-3}{2}^{-1}[10\C{t}{5} -  \varphi s]$},
$$
which, combined with (\ref{bigB}), gives
$
s< 10\cone\varphi^{-1}\Cc{t}{5}.
$\qqed

\medskip
We may then finish {\em via} the following simple lemma.
Let $k,l$ be integers with $k<l$ and
$W$ a set of size $t$.  Suppose that
for each $R\in \C{W}{k}$
we are given some
$\gs_R:R\ra\{0,1\}$, and for $R,S\in \C{W}{k}$
write $R\sim S$ if $\gs_R$ and $\gs_S$ agree on $R\cap S$.
Say $L\in \C{W}{l}$ is {\em consistent} if $R\sim S$
$\forall\,~R,S\in \C{L}{k}$.
\begin{lemma}\label{Lproper'}
For all k, l as above and $\eps>0$ there is a $\xi>0$ such that
(with notation as above) if at least $(1-\xi)\C{t}{l}$
$l$-subsets of $W$ are consistent, then there
is some $f:W\ra \{0,1\}$ such that
$\gs_R\equiv f|_R $ for all but at most $\eps\C{t}{k}$
$k$-subsets $R$ of $W$.
\end{lemma}
\nin
We will prove this only for $k=3$ and $l=5$, in which case
we may take $\xi =(\eps/6)^2$.
The proof of the general case, an induction on $k$, is in a similar
vein, though not exactly a generalization of the argument given
here.

\medskip
Of course to get Lemma \ref{step3}(b) from (the case $k=3,l=5$ of)
Lemma \ref{Lproper'}
we take $W$ to be the set of blocks of $\pee^*$,
set $\gs_R=\pi_P$ whenever $P$ is a proper triad and $R$ its
set of blocks, and define
$\gs_R$ arbitrarily for the remaining $R$'s.
(Here we use \eqref{c1phi}.)\qqqed

\mn
{\em Proof of Lemma} \ref{Lproper'} (for $k=3,l=5$).
Let $\xi$ be as above,
set $\ga = \frac{1}{8}\sqrt{\xi}$, and
say $x\in W$ is {\em bad} if there are at least $\ga (t)_4$
pairs $\{R,S\}$ with: $R,S\in \C{W}{3}$;
$R\cap S=\{x\}$;
and
$R\not\sim S$.
If the number of bad $x$'s is $b$ then the number of inconsistent
$5$-sets is at least $\frac{1}{15}b\ga (t)_4$,
so we have
$b< \frac{15}{\ga (t)_4}\xi \C{t}{5} <\frac{\xi}{8\ga}t$.

If, on the other hand, $x$ is not bad then (by Lemma \ref{component})
there is $f(x)\in \{0,1\}$ such that $\gs_R(x)=f(x)$
for at least (say) $(1-8\ga)\C{t-1}{2}$
$3$-sets $R\ni x$.
So extending this $f$ arbitrarily to the bad $x$'s
we find that the number of 3-sets $R$ that fail to satisfy
$\gs_R\equiv f|_R$ is at most
$t\cdot 8\ga\C{t-1}{2} +b\C{t-1}{2} <
(8\ga + \frac{\xi}{8\ga})t\C{t-1}{2} = \eps\C{t}{3}$.\qed

\mn
{\em Proof of Lemma} \ref{ML1}.
We first show that clauses not belonging to $\Cl(\pee^*)$
are more or less irrelevant.
We are interested in the number of possibilities for
$\cee\sm\Cl(\pee^*)$ with $\cee\sim \pee^*$.
Members of $\cee\sm \Clp$ are either

\mn
(i) clauses not supported on triads of $\pee^*$ or

\mn
(ii)  clauses belonging to patterns $\pi$ that are supported
on triads of $\pee^*$, but that are not patterns of $\pee^*$
(i.e. for which $\ddd_{\pi}\leq 2d_0$).

\mn
The total number of {\em possible} clauses
of the first type is
$O(\gd + \eps_1+t^{-1})n^3=O(\gd n^3)$
(see \eqref{gdt}),
where the first term, given by
\eqref{RegPart'}, is for clauses
supported on triads of the underlying partition $\pee$ that are
not triads of $\pee^*$.
(The other two terms bound the number of clauses that use either $V_0$
or some $P^{ij}_0$, or that meet some block more than once.)
On the other hand, no
$\cee\sim \pee^*$ contains more than $16d_0\C{n}{3}$ clauses of type (ii).
Thus we have (using $\sum_{i\leq k}\C{m}{i}\leq \exp[H(k/m)m]$)
\beq{NvsN*}
N^*(\pee^*) < \exp[8H(2\ddd_0)\Cc{n}{3}+ O(\gd )n^3]N(\pee^*).
\enq
Thus (\ref{bigP})
implies
\begin{eqnarray}\label{bigP'}
h(\pee^*) 
&> &[(1 + 5\eps_2l)m^3l^{-3}]^{-1}[(1-\ctwo)\Cc{n}{3}
-8H(2\ddd_0)\Cc{n}{3}- O(\gd )n^3]\nonumber\\
&> &(1-2\ctwo)\Cc{t}{3}l^3
\end{eqnarray}
(where we used
$\ctwo >> \max\{H(2\ddd_0), \gd,\eps_2l\} $ ($=H(2\ddd_0)$)
and $\C{n}{3}> \C{t}{3}m^3$).

\medskip
We next observe that \eqref{bigP'}
(and so (\ref{bigP})) implies
\beq{mostBCs}
\mbox{all but at most $2\ctwo\cone^{-1}l^{\C{t}{2}}$
BC's of $\pee^*$ are coherent.}
\enq
{\em Proof.}  This is similar to the proof of
\eqref{beta5}.
By Lemma \ref{step3}(b), the number of
incoherent BC's of $\pee^*$ is at most
$$
s:=|\{\gb:\mbox{$\gb$ a BC of $\pee^*$;
$h(\gb)< (1-\cone)\C{t}{3}$}\}|.
$$
Thus Corollary \ref{Shcor}(b) and Lemma \ref{step3}(a) give
\begin{eqnarray*}
h(\pee^*) &\leq &
l^{-\C{t}{2}+3}\sum \{h(\gb):\mbox{$\gb$ a BC of $\pee^*$}\}\\
&<& l^{-\C{t}{2}+3}( (l^{\C{t}{2}}-s)\Cc{t}{3}+s(1-\cone)\Cc{t}{3}),
\end{eqnarray*}
which with \eqref{bigP'} implies
$
s< 2\ctwo\cone^{-1}l^{\C{t}{2}}.
$\qqed

\medskip
For the rest of this argument $\gb$ ranges over
BC's (of $\pee^*$), $P$ and $Q$ over
{\em triads of} $\pee$, and $A,B,C$ over blocks.
For each coherent $\gb$ we fix some $f_{\gb}$ as in \eqref{coherence'}
and assign an arbitrary (convenient but irrelevant)
$f_{\gb}:\{\mbox{blocks}\}\ra \{0,1\}$ to each
incoherent $\gb$.

Say $P$ and $Q$ {\em disagree} at a common block
$A$ if at least one of $P,Q$ is not proper or
(both are proper and)
$\pi_{_P}(A)\neq\pi_{_Q}(A)$.
(Here one should think of $P$ and $Q$ as having just the one block in common;
effects due to pairs with larger overlap will be insignificant.)
We now proceed roughly as follows.  An averaging argument
shows
that for most blocks $A$ there are few pairs $P,Q$ that disagree at $A$.
When this happens there must be a value for $f(A)$ that
agrees
with most of the triads using $A$.
The remaining few $f$-values are then of no concern and may be assigned
arbitrarily.

To say this properly, write $P\not\sim _A Q$ if $P$ and $Q$
disagree at $A$ and have no other block in common.
Write
$P\not\sim_A\gb$ if $P$ is a triad of $\gb$ and either
$P$ is improper
or $\pi_P$ disagrees with $f_{\gb}$ at the
block $A$ of $P$, and $P\not\sim \gb$ if $P\not\sim_A\gb$
for some block $A$ of $P$.
Setting
$$M =|\{(\gb,P,Q,A):\mbox{$P,Q$ triads of $\gb$; $P\nsa\gb$ or $Q\nsa\gb$}\}|,$$
 we have
\begin{eqnarray}
M
&\leq &
2\mbox{$\C{t-1}{2}|\{(\gb,P,A):P\nsa\gb\}| $}\nonumber\\
&\leq &
\mbox{$6\C{t-1}{2}|\{(\gb,P):P\not\sim \gb\}|$}\nonumber \\
&\leq& \mbox{$6\C{t-1}{2}(2\gz_1 )\lct\ct $}
<  \mbox{$O(\gz_1 t^5 \lct)$}\nonumber,
\end{eqnarray}
where we use $\gz_1$ to bound both the fraction of
incoherent $\gb$'s (see \eqref{mostBCs} and \eqref{gz1c2c1})
and the fraction of triads
that disagree with $f_{\gb}$ when $\gb$ is coherent.
But we also have
$$
M\geq |\{(A,P,Q):P\nsa Q\}|l^{\C{t}{2}-6};
$$
thus
$$
\sum_A|\{(P,Q):P\nsa Q\}| =
|\{(A,P,Q):P\nsa Q\}| < O(\gz_1 t^5l^6),
$$
implying
\beq{goodA}
|\{(P,Q):P\nsa Q\}| <\sqrt{\gz_1}t^4l^6
\enq
for all but at most $O(\sqrt{\gz_1} t)$ $A$'s.

For $A$ satisfying (\ref{goodA})
we again appeal to Lemma \ref{component},
applied to the graph $G=G_A$ having vertices the triads (of $\pee$)
that use $A$, and $PQ$ an edge if $P,Q$ are proper and
$\pi_P(A)=\pi_Q(A)$
(so improper triads become isolated vertices).
We have
$|V(G)|=\C{t}{2}l^3$ and
$|E(\ov{G})| <\sqrt{\gz_1}t^4l^6 + t^3l^6 $
(the negligible second term being a bound on the number of pairs $P,Q$
that share at least one additional block);
so the lemma says there
is
some $f(A)\in \{0,1\}$ such that
$\pi_{_P}(A)=f(A) $
for all but at most $O(\sqrt{\gz_1} t^2l^3)$ triads $P$ using $A$.

Finally, extending this $f$ arbitrarily to $A$'s failing (\ref{goodA}),
we find that the number of triads (of $\pee$) that are improper or
disagree with $f$---so in particular the number (needed for
\eqref{coherence})
that {\em are}
proper and disagree with $f$---is less than
$O(\sqrt{\gz_1}t^3l^3)$; so, in view of \eqref{gz1gz2},
$\pee^*$ is coherent.\qed

\section{Proof of Lemma \ref{Pw}}\label{Witnesses}

It will now be convenient to work with triangles rather than triads,
which we can arrange, e.g., by observing that \eqref{coherence}
implies
\beq{coherence''}
\mbox{{\em all but at most $2\gz_2 \C{t}{3}m^3$ triangles belong
to triads that
agree with $f$}}
\enq
(by \eqref{RegPart'},
since $\gd$ is much smaller than $\gz_2$).

\mn
We first need to show that
$f$ as in \eqref{coherence''}
is mostly 1.
Say (just for the present argument) that a block $V_i$ is ``bad" if
at least $.05 \C{t-1}{2}m^3$ triangles belong to triads that disagree
with $f$ at $V_i$.
Let $M$ be the number of bad $V_i$'s and
$N$ the number of pairs $(V_i,K)$ with
$V_i$ a block of $\pee^*$ and $K$
a triangle belonging to a triad
that disagrees with $f$ at $V_i$.
Then
$$ 6\gz_2 \Cc{t}{3}m^3 \geq N \geq .05 M  \Cc{t-1}{2}m^3$$
gives
$M\leq 40\gz_2t$.

Suppose, on the other hand, that $V_i$ is good (i.e. not bad).
Then
the number of clauses (of $\cee$) that agree with $f$
at $V_i$ is at least
$\frac{1}{3}(.95) \C{t-1}{2}m^3$
(since each triad $P$ that agrees with $f$ at $V_i$
is proper and thus
contributes at least $\frac{1}{3}t(P)$ such clauses),
while the number that disagree
is at most $4(.05 +\ddd_0)\C{t-1}{2}m^3$.
There is thus (since $\frac{1}{3}(.95) > 4(.05+\ddd_0)$)
some $x\in V_i$ that belongs to more clauses
that agree with $f$ at $x$ than that disagree, so that $m(x)\geq m(\bar{x})$
implies that $f(V_i)=1$.
So we have shown that
$$|f^{-1}(0)|\leq 40\gz_2 t.$$

Now suppose for a contradiction that
$\ww$ is a witness for some $C\in \cee$
and $|\ww^{-1}(1)|>\gz n$.
Then for the set, say $\W$, of blocks $V_i$
satisfying
$$
\mbox{$f(V_i)=1$ and
$|\ww^{-1}(1)\cap V_i|> \gz m/2$,}
$$
we have
$$ \gz n ~<~ |\ww^{-1}(1)|~<~40\gz_2  n +|\W| m + \gz n/2,$$
whence
$$|\W|~\geq  ~(\gz /2 -40\gz_2 )n/m ~\geq~(\gz /2 -40\gz_2 )t.$$
It then follows from \eqref{coherence''}, using (say)
\beq{param1}
(\gz /2 -40\gz_2)^3 >3\gz_2,
\enq
that there is some
triad $P$ that agrees with $f$, all three of whose blocks
are in $\W$ (which, note, implies $\pi_{_P} \equiv 1$).
But then
$$
(1-8\eps_2 l)(\gz/2)^3 > \gd (1+5\eps_2l^3)
$$
(implied by \eqref{param1}) and
$(\gd,r)$-regularity of $P$
imply that there is some $C\neq xyz\in \cee$ supported by $P$,
so that
$\ww$ cannot have been a witness.
(In more detail:  Suppose the blocks of $P$ are $V_i,V_j,V_k$,
and let $V_u'= \ww^{-1}(1)\cap V_u$.
Then using
Proposition \ref{Ptriangles'} (both the upper and lower bounds),
we find that for the subtriad $Q$ of $P$ spanned (in the obvious sense)
by $V_i',V_j',V_k'$, we have
$$
|T(Q)|> (1-8\eps_2 l)(\gz/2)^3m^3l^{-3} > \gd (l^{-3}+5\eps_2)m^3>\gd t(P);
$$
thus
$(\gd,r)$-regularity (here $r=1$ would suffice) gives
$\ddd_{\pi_P}(Q) > \ddd_{\pi_P}-\gd$, implying the existence of $xyz$ as
above.\qed

\section{Recursion}\label{Phase2}

\mn
Here we prove (\ref{I*}).
From this point we write simply $X$ for $X_n$ (the set of variables),
and use $a,b,c,u,v,w,x,y,z$ for members of $X$.
\glossary{name={positive clause},description={Contains only positive literals.}}
\glossary{name={negative clause},description={Contains only negative literals.}}
We call a clause {\em positive (negative)} if it contains
only positive (negative) literals, and
{\em non-positive} if it contains at least one negative
literal.
We assume throughout that all $\cee$'s under discussion belong to $\ii$
(and, as usual, that $n$ is large enough to support our assertions).

As the form of (\ref{I*}) suggests, the proof will proceed by
removing from $\ii$ $\cee$'s exhibiting various ``pathologies,"
eventually leaving only (a subset of all)
$\cee$'s containing only positive clauses;
these account for the main term,
$2^{\C{n}{3}}$, on the right hand side of (\ref{I*}).

The arguments again involve interplay of a number of small constants,
and we begin by naming these and specifying what we will assume
in the way of relations
between them.  In addition to $\cc$ (from (\ref{I*})) and $\gz$
(from \eqref{zeta}), we will use constants $\ga$, $\vt$ and
$\xi$, assumed to satisfy the (satisfiable) relations
\beq{cbds}
0<\cc < \min\{\xi, \vt^3-7H(2\gz),
\tfrac{2-\log 3}{12} - 3H(\vr/3)
\}=\vt^3-7H(2\gz),
\enq
where $\vr = \sqrt{2\ga}+\gz$,
and
\begin{eqnarray}\label{xibd}
\xi &< &\min\{\ga -2\vt,\sqrt{.04-2\vt}-\vartheta,
0.1-7 H(2\gz),\nonumber\\
&& ~~~~~~~~
1-\tfrac{1}{3}H(\tfrac{1}{10}) - 0.3 \log 7-7H(2\gz+\ga)\}
=\ga -2\vt.
\end{eqnarray}
(These hold if all parameters are small and, for example,
$\ga > 2\xi>5H(\vartheta)$ and $\vartheta > 7H(2\gz)$.)

\mn
{\em Step} 0.
Let
$$\ii_1= \{\cee\in \ii:
\mbox{each variable is used at least $\frac{1}{10} \C{n-1}{2}$ times in $\cee$}\}.
$$
Then
\beq{II1}
\mbox{$|\ii\sm \ii_1|< \exp[.8\C{n}{2}]I(n-1)$.}
\enq

\mn
{\em Proof.}
There are at most

$$\mbox{$n \sum\{\C{8\C{n-1}{2}}{t}:t\leq \frac{1}{10}\C{n-1}{2}\} ~<~
\exp[H(\frac{1}{80})8\C{n-1}{2}] ~<~
\exp[.8\C{n}{2}]$}$$

\mn
ways to choose a variable $x$ to be used fewer than
$\frac{1}{10} \C{n-1}{2}$ times, together with the clauses that use $x$,
and the collection
of clauses of $\cee$ not using $x$ is an (irredundant)
formula on the $n-1$ remaining variables.\qqed

\mn
{\em Step} 1.
If $\cee\in \ii$ then for any two variables $u,v$
there are at most $\gz n$ variables $w$ for which
$uv\bar{w}\in \cee$.
The same bound applies to $w$'s with
$u\bar{v}\bar{w}\in\cee$
and those with $\bar{u}\bar{v}\bar{w}\in\cee$.

\mn
{\em Proof.}
If $\ww$ is a witness for $uv\bar{w}\in\cee$ then
any $x\neq w$ with $uv\bar{x}\in\cee$ must lie in
$\ww^{-1}(1)$.
The other cases are similar.\qqed

\mn
In particular:

\mn
(a)  for any $u$, $\cee$ contains at most $\gz n^2$ clauses
of each of the forms $uv\bar{w}$, $u\bar{v}\bar{w}$,
$\bar{u}v\bar{w}$, $\bar{u}\bar{v}\bar{w}$;

\mn
(b)  $\cee$ contains at most (say) $2\gz n^3$ non-positive
clauses;

\mn
(c)  if $\cee\in \ii_1$
then, for any $u$, $\cee$ contains at least (say) $0.02 n^2$
positive clauses using $u$ (by (a), since $\cee\in \ii_1$ implies
$m(u)\geq \frac{1}{20} \C{n-1}{2}$).

\bn
{\em Step} 2.
Let
$\ii_2$ consist of those $\cee\in \ii_1$ that satisfy
\beq{ubar}
\mbox{for each $u$,
$\cee$ contains at most $\ga n^2$ clauses $\bar{u}vw$.}
\enq
Then
\beq{I1I2}
\mbox{$|\ii_1\sm \ii_2|< \exp[(1-c)\C{n}{3}]
+ \exp[(1-c)\C{n}{2}]I(n-1)$.}
\enq
{\em Proof.}
We should show that the number of $\cee$'s in $\ii_1$ violating
(\ref{ubar}) is at most the right hand side of (\ref{I1I2}).
Given such a $\cee$ we
fix $u$ violating (\ref{ubar}) and set $Y= X\sm \{u\}$,
\[
\mbox{$R=\{\{a,b\}\sub Y:uab\in \cee\}$,
$~B=\{\{a,b\}\sub Y:\bar{u}ab\in \cee\}$,}
\]
\[
\mbox{
$S=\{a\in Y: d_R(a)\leq \vartheta n\}$, $~T=\{a\in Y: d_B(a)\leq \vartheta n\}$}
\]
(where
we regard $R$ and $B$ as graphs on $Y$ and use
$d$ for degree) and
$Z=Y\sm (S\cup T)$.

The main point here is
that, because $\cee$ is irredundant,
\beq{umain}
\mbox{if $ab\in R$ and $ac\in B$ (and $b\neq c$) then $abc\not\in \cee$.}
\enq

Since the number of clauses $\bar{u}vw$,
which we are assuming to be at least $\ga n^2$,
is at most $(n-|T|)n +|T|\vartheta n \leq (|S|+|Z|)n +\vartheta n^2$,
we must have {\em either} $|Z|> \vartheta n$ {\em or}
$|Z|\leq \vartheta n$ and $|S|> \xi n$
(see \eqref{xibd}).

Suppose first that
$|Z|>\vartheta n$.
In this case, once we have specified $Z$ and the $R$-
and $B$-edges meeting $Z$, (\ref{umain}) gives at least
$\vartheta n\cdot \vartheta n\cdot (\vartheta n-1)/6$ positive clauses
$abc$ that are known to {\em not} belong to $\cee$.
We may thus (crudely) bound the number of possibilities for $\cee$
of this type by the product of the factors:
$n$ (corresponding to the choice of $u$);
$2^n$ (choose $Z$);
$\exp[n^2]$ (for the $R$- and $B$-edges meeting $Z$);
$\exp[H(2\gz) \cdot7\C{n}{3})]$ (for the remaining
non-positive members of $\cee$ (i.e. those not of the
form $\bar{u}vw$); here we use (b) of Step 1);
and $\exp[(1-\vartheta^3)\C{n}{3}]$
(for the remaining positive members of $\cee$).
This product is less than the first term on the right hand
side of (\ref{I1I2}).

Next suppose
$|Z|\leq\vartheta n$ and $|S| >\xi n$.
We first observe that
$n-|S|$ can't be
{\em too} small:
the number of positive clauses of $\cee$ using $u$
is at least $0.02 n^2$ (by (c) of Step 1), but also
at most
$|S|\vartheta n + \C{n-|S|}{2}$, which, after a little
calculation, gives
$n-|S|> \sqrt{.04-2\vartheta}~n$.
Thus in the present case we must have
$|T|> (\sqrt{.04-2\vartheta}-\vartheta)n>\xi n$.

We may specify a $\cee$ of the present type
(i.e. with $|Z|\leq \vartheta n$ and $|S|> \xi n$,
so also $|T|>\xi n$)
by choosing:
(i)  $u$;
(ii)  $S$ and $T$ (so also $Z$);
(iii)  the $R$-edges meeting $S\cup Z$ and
the $B$-edges meeting $T\cup Z$;
(iv)  the $R$-edges contained in $T':=T\sm S$
and
the $B$-edges contained in $ S':=S\sm T$; and
(v) the clauses not involving the variable $u$.
The numbers of choices in (i), (ii) and (v) are
at most $n$, $4^n$ and $I(n-1)$ (respectively), while those
for for (iii) and (iv) are bounded by
$$
\mbox{$\exp[2\vartheta n^2 +(|S|+|T|)H(\vartheta) n +
\C{|S'|}{2}+\C{|T'|}{2}]$}.
$$
Combining these bounds with the easy
$$
\mbox{$\C{|S'|}{2}+\C{|T'|}{2} < \C{n-1}{2} -\xi(1-\xi)n^2
$,}
$$
we find that the number of $\cee$'s in question
is less than
$$
n4^n\exp[\Cc{n}{2} -(\xi(1-\xi) -2\vartheta -2H(\vartheta))n^2]I(n-1),
$$
which is less than the second
term on the right hand
side of (\ref{I1I2}).\qqed

\mn
Note that $\cee \in \ii_2$ implies (by (a) of Step 1) that for any $u$,
\beq{nonpos}
\mbox{$\cee$ contains at most $(4\gz +\ga)n^2$
non-positive clauses using $u$ or $\bar{u}$.}
\enq

\bn
{\em Step} 3.
For a
variable $u$,
set
$\XX_u = \{\{v,w\}:uvw\in \cee\}$
and $\bar{\XX}_u =\C{X\sm \{u\}}{2}\sm \XX_u$.
Let
$\ii_3$ consist of those $\cee\in \ii_2$
with the property that
for any three variables $u,v,w$,
$$
\mbox{each of $
|\XX_u\cap\XX_v\cap \XX_w|,~
|\XX_u\cap\XX_v\cap \bar{\XX}_w|,~
|\XX_u\cap \bar{\XX}_v\cap \bar{\XX}_w|$} ~~~~~~~
$$
\beq{X's}
~~~~~~~~\mbox{and $|\bar{\XX}_u\cap \bar{\XX}_v\cap \bar{\XX}_w|$
is at least
$0.1 \C{n}{2}$.}
\enq
(The ``$0.1$" is just a convenient constant smaller than $1/8$.)
We assert that
\beq{I2I3}
\mbox{$|\ii_2\sm \ii_3|< \exp[(1-\cc)3\C{n}{2}]I(n-3)$.}
\enq
{\em Proof}.
We may choose $\cee\in \ii_2\sm \ii_3$ by choosing:

\mn
(i) $u,v,w$ violating (\ref{X's});

\mn
(ii)  the non-positive clauses involving at least one of $u,v,w$;

\mn
(iii)  the positive clauses involving $u,v,w$;

\mn
(iv)  the clauses not involving $u,v,w$.

\medskip
The numbers of possibilities for the choices in
(i), (ii) and (iv) may be bounded by $\C{n}{3}$,
$\exp[3H((4\gz +\ga)n^2)/(7\C{n}{2})\cdot 7\C{n}{2}]<
\exp[21H(2\gz +\ga)\C{n}{2}]$ (see (\ref{nonpos}))
and $I(n-3)$ respectively.
The main point is the bound for the number of choices in (iii),
which,
apart from the $2^{O(n)}$ possibilities for clauses involving at
least two of $u,v,w$, is bounded by the number of choices
for an ordered
partition of $\C{X\sm \{u,v,w\}}{2}$ into eight parts,
at least one of which has size less than $0.1\C{n}{2}$.
We assert (a presumably standard observation) that this number is
less than
$
8\exp[(H(.1) + .9\log 7 )\Cc{n}{2}].
$
which finishes Step 3 since the product of the preceding bounds
is less than the right hand side of (\ref{I2I3}).

For the assertion, notice that
the log of the number of (ordered) partitions
$[m]= Z_1\cup \cdots \cup Z_8$ with $|Z_1|< 0.1 m$ is
$H(Y_1\dots Y_m)\leq \sum H(Y_i)$,
where we choose $(\ZZ_1\dots \ZZ_8)$ uniformly from the set of
such partitions and set $Y_i=j$ if $i\in \ZZ_j$.
(The inequality, an instance of Lemma \ref{Shearer},
is a basic (easy) property of entropy; see e.g. \cite[Theorem 2.6.6]{CT}.)
Setting $p_i(j)=\Pr(Y_i=j)$ ($=\Pr(i\in \ZZ_j)$) and
$\bar{p}_j =m^{-1}\sum_ip_i(j)$, we have
\begin{eqnarray*}
\mbox{$\sum H(Y_i)$} &=& \mbox{$\sum_j\sum_ip_i(j)\log \tfrac{1}{p_i(j)}$}\\
&\leq &
\mbox{$m\sum_j \bar{p}_j\log \tfrac{1}{\bar{p}_j}  =
m H(\bar{p}_1\dots \bar{p}_8)$}
\end{eqnarray*}
(by
Jensen's Inequality) and
$$
H(\bar{p}_1\dots \bar{p}_8)
\leq H(\bar{p}_1) + (1-\bar{p}_1)\log 7 <
H(0.1) + 0.9\log 7$$
(using $H(X)\leq \log |{\rm range}(X)|$ for the first inequality).
\qqed

\mn
{\em Step} 4.
Let
$$\ii_4=\{\cee\in \ii_3:
\mbox{no clause of $\cee$ uses more than one negative literal.}\}
$$
Then
\beq{I3I4}
\mbox{$|\ii_3\sm \ii_4|< \exp[(1-c)\C{n}{3}]$.}
\enq
{\em Proof.}
We first observe that $\cee\in \ii_3$
{\em cannot} contain a clause with exactly two
negative literals.
For suppose $\bar{u}\bar{v}w\in\cee$.
Since $\cee\in \ii_3$, there is some pair $\{a,b\}$
with $abu,abv,abw\in\cee$;
but this is impossible, since
a witness for $abw$ must agree with
at least one of $abu$, $abv$, $\bar{u}\bar{v}w$.

While the preceding argument doesn't quite work to exclude
negative
clauses,
the assumption that $\bar{u}\bar{v}\bar{w}\in\cee$ is extremely
restrictive, since it says that whenever
$\{a,b\}\in \XX_u\cap\XX_v\cap \XX_w$, there cannot be {\em any}
$c\not\in \{u,v,w\}$ with
$abc\in \cee$
(since a witness for $abc$ would have to agree with one of
$abu$, $abv$, $abw$, $\bar{u}\bar{v}\bar{w}$).
So we may bound the number of $\cee$'s that do contain negative
clauses by the product of:
$n^3$ (choose $u,v,w$); $\exp[n^2]$ (choose $\XX_u\cap\XX_v\cap \XX_w$);
$\exp[7H(2\gz)\C{ n}{3} + O(n^2)]$ (for clauses that either are non-positive or
involve $u,v$ or $w$; here we again use (b)); and
$\exp[ \C{n-3}{3} - 0.1\C{n}{2}(n-3)/3]< \exp [.9\C{n}{3}]$
(for the remaining positive clauses;
here the subtracted term corresponds to
triples known to contain members of $\XX_u\cap\XX_v\cap \XX_w$).
And again, the product of these bounds is less than
$\exp[(1-c)\C{n}{3}]$.\qqed

\mn
{\em Step} 5.
Finally, we set
$$\ii_5=\{\cee\in \ii_4:
\mbox{$\cee$ contains no clause with exactly one negative literal}\}
$$
(so $\ii_5 \sub \{\cee\in \ii:\mbox{$\cee$ contains
only positive clauses}\}$)
and show
\beq{I4I5}
\mbox{$|\ii_4\sm \ii_5|< \exp[\C{n}{3}-cn]$.}
\enq

\mn
{\em Proof.}
We show that for any $\ttt >0$
(by (b) of Step 1 $\ttt$ will be at most
$\gz n^3$, but we don't need this),
\beq{I5}
|\{\cee\in \ii_4:\mbox{$\cee$
has exactly $\ttt$ non-positive clauses}\}|
<\mbox{$\exp[\C{n}{3}-c'n]$}
\enq
for a suitable $c'$; this gives (\ref{I4I5})
for any $c<c'$.

Fix $\ttt$ and suppose $\cee$ is as in (\ref{I5}).
The main point driving the argument (which, however,
will take us a while to get to)
is:
\beq{Cmp}
\mbox{if $\bar{u}vw\in \cee$ and $a\not\in \{u,v,w\}$,
then $|\cee\cap \{auv,avw\}|\leq 1$}
\enq
(since a witness for $avw$ must agree with either $auv$ or $\bar{u}vw$).

Let $\cee'$ be the set of non-positive clauses in $\cee$.
It will be helpful to introduce an auxiliary collection:
for each $C\in \cee'$, we will fix an ordering of the three literals
in $C$ with the negative literal first,
and write $\cee''$ for the resulting collection of ordered triples.
We assert that we can do this
so that
\beq{few}
|\{w:(\bar{u},v,w)\in \cee''\}|\leq  \sqrt{\ga/2}~n ~~~~~\forall\,u,v.
\enq
\sugg{
Let $\cee'$ be the set of non-positive clauses in $\cee$.
Here it will be helpful to slightly change our notation:
for each $C\in \cee'$, we will fix an ordering of the three literals
in $C$, with the negative literal first
and
regard $\cee'$ as a set of
{\em ordered} clauses (so
``$\bar{u}vw\in \cee'$" is no longer the same as
``$\bar{u}vw\in \cee'$").
We assert that we do this
so that
\beq{few}
|\{w:\bar{u}vw\in \cee'\}|\leq  \sqrt{\ga/2}~n ~~~~~\forall\,u,v.
\enq
}

This will follow from
\mn
\begin{prop}\label{Hatami}
Any (simple) graph admits an orientation with all out-degrees
at most $\sqrt{|E(G)|/2}$.
\end{prop}

\mn
{\em Proof (sketch).}
A precise statement
(due to Hakimi \cite{Hakimi}; see also
\cite[Theorem 61.1, Corollary 61.1b]{Schrijver})
is:
for any graph $G =(V,E)$ and
$c:V \ra \Nn$,
there is an orientation with
$d^+_v\leq c_v ~\forall\,v$
(where, of course, $d_v^-$ is
the out-degree of $v$)
iff for every $W\sub V$,
$
|E(G[W])| \leq \sum\{c_v:v\in W\}
$;
in particular, there is an orientation with
$d^+_v\leq c ~\forall\,v$ iff
$c\geq \max\{|E(G[W])|/|W|:W\sub V\}$,
which is easily seen to hold with
$c = \lceil \sqrt{|E(G)|/2}~\rceil$.

(Alternatively
it's easy to see that
orienting each edge toward the end of larger
degree (breaking ties arbitrarily) gives maximum out-degree less
than $\sqrt{2|E(G)|}$, which would also be fine for present purposes.)\qqed

\medskip
To get \eqref{few} from Proposition \ref{Hatami}, regard, for a given $u$,
$\{vw:\bar{u}vw\in \cee'\}$ as the edge set of
a graph $G_u$ on $X\sm \{u\}$, and choose an orientation of
$E(G_u)$ as in the proposition.
We have
$|E(G_u)|\leq \ga n^2$ (by (\ref{ubar})); so
interpreting orientation of $vw$ toward $w$ as specifying
$(\bar{u},v,w)\in \cee'$ gives \eqref{few}.\qqed

\medskip
Of course there will typically be many choices of $\cee''$ as above,
and we fix one such for each $\cee'$.
Given $\cee''$,
set $\g =\g(\cee') =\{\{\{u,v\},\{v,w\}\}:(\bar{u},v,w)\in \cee''\}$.
Regard $\g$ as a multigraph on the vertex set $\C{X}{2}$,
and let $\nu$ and $\tau$ denote its matching and (vertex) cover
numbers.  Then
\beq{taunu}
2\nu\geq \tau \geq \lceil\frac{\ttt}{\vr n}\rceil
\enq
(where, recall, $\vr =\sqrt{2\ga} +\gz$).
Here the first inequality is standard (and trivial)
and the second follows from the fact that $\g$ has $\ttt$ edges
and maximum degree at most $\vr n$, the latter by
(\ref{few}) and Step 1.

We now consider the number of possibilities for $\cee$
with given a $\ttt $, $\tau$ and
$\nu$.
We first specify $\cee'$ by choosing a vertex cover $\T$ for the
associated $\g$ and then a collection of $\ttt$ clauses,
each using (the variables from) at least one member of $\T$.
The number of possibilities for these choices is at most
$\C{n^2}{\tau}\C{3\tau n}{\ttt}$.

We now suppose $\cee'$ has been determined and consider
possibilities for the set, say $\cee_0$ ($=\cee\sm\cee'$),
of positive clauses of $\cee$.
Let $\m$ be some maximum matching of $\g$,
say $\m=\{\{\{u_i,v_i\},\{v_i,w_i\}\}:i\in [\nu]\}$.
(We could specify $\bar{u}_iv_iw_i\in \cee'$, but this is now
unnecessary.)

Let $\J$ be the set of all pairs of 3-sets
$\{\{a,u_i,v_i\},\{a,v_i,w_i\}\}$
such that $\{\{u_i,v_i\},\{v_i,w_i\}\}\in \m$ and $a\not\in \{u_i,v_i,w_i\}$,
and let $\K$ be the set of 3-sets belonging to pairs in $\J$.
Then $\J$ is a set of at least $\nu(n-3)/2$ pairs of 3-sets
(a given pair $\{\{x,y,z\},\{x,y,w\}\}$ can arise
with $x$ in the role of $v_i$ and $y$ in the role of $a$
or vice versa)
with the property that
no 3-set belongs to more than three members of $\J$
(since $\m$ is a matching);
so in particular $|\K|\geq \nu(n-3)/3$.

We assert that the number of possibilities for
$\cee_0\cap \K$ is at most
$$
\mbox{$\exp[\C{n}{3}-\frac{1}{6}\nu(n-3)(2-\log 3)].$}
$$

\mn
{\em Proof.}
This is another (somewhat more interesting)
application of
Lemma \ref{Shearer}.
Let $W=\C{X}{3}$
(thought of as the collection of possible positive clauses);
let $\f$ be the collection of possible $\cee_0$'s
(compatible with the given $\cee'$);
and let $\h$ consist of all pairs from $\J$
(note these are now pairs of elements of $W$) together with,
for each $T\in W$, $3-\eta(T)$ copies of the singleton $\{T\}$,
where $\eta(T)\leq 3$ is the number of times $T$ appears as a
member of some pair in $\J$.
As noted earlier the key point is (\ref{Cmp}), which in the
present language says that no member of
$\f$ contains any $\{S,T\}\in \J$.
This implies in particular that for each such $\{S,T\}$, we have
$|{\rm Tr}(\f,\{S,T\})|\leq 3$,
so that Lemma \ref{Shearer} gives
\begin{eqnarray*}
\log |\f| &\leq & \mbox{$\frac{1}{3}\left[
\sum_{T\in W}(3-\eta(T)) + |\J|\log 3\right]$}\\
&\leq& \mbox{$\C{n}{3} - \frac{1}{6}\nu (n-3)(2-\log 3)$}
\end{eqnarray*}
(since $\sum \eta(T) =2|\J|$ and $|\J|\geq \nu(n-3)/2$).\qqed

\medskip
Finishing the proof of (\ref{I4I5}) is now easy.
We have shown that the number of
possibilities for $\cee$
with given $\ttt $, $\tau$ and
$\nu$ is at most
\[
\mbox{$\C{n^2}{\tau}\C{3\tau n}{\ttt}
\exp[\C{n}{3}-\frac{1}{6}\nu(n-3)(2-\log 3)]$}
~~~~~~~~~~~~~~~~~~~~~~~~~~~
\]
\[
~~~~~~~~~~~
\mbox{$< \exp\left[\C{n}{3} +\left\{\log\frac{en^2}{\tau}
+3nH(\vr/3)-
\frac{(n-3)(2-\log 3)}{12}\right\}\tau\right]$}
\]
(where we used (\ref{taunu}) (second and first inequalities
respectively) for the last two terms in
the exponent), and
summing over $\tau$ and $\nu$
shows that the left side of (\ref{I5}) is less than
$\exp[\C{n}{3}- c'n]$ for any
$c' < (2-\log 3)/12  -
3H(\vr/3)$.\qqed

\medskip
Finally, combining (\ref{II1}), (\ref{I1I2}), (\ref{I2I3}),
(\ref{I3I4}) and (\ref{I4I5}) (and, of course, the fact that
$|\ii_5|\leq \exp[\C{n}{3}]$)
gives (\ref{I*})
(where we again absorb terms $\exp [(1-c)\C{n}{3}]$
from (\ref{I1I2}) and (\ref{I3I4})
in the term $\exp[\C{n}{3}-\cc n]$).\qed

\bn
$~$

\bn
{\bf Acknowledgments}
Thanks to Dan Kleitman for some stimulating
conversations and
to Nick Wormald for references \cite{Hakimi} and \cite{Schrijver}.
Part of this work was carried out while the second author was
visiting MIT.

\bigskip
\begin{flushleft}
Department of Mathematics\\
Rutgers University\\
Piscataway NJ 08854 USA\\
ilinca@math.rutgers.edu\\
jkahn@math.rutgers.edu
\end{flushleft}

\end{document}